\documentclass[10pt]{article}
\usepackage{amsmath, amsthm, amssymb, amsfonts, mathrsfs, mathtools, upgreek}
\usepackage{graphicx}
\usepackage{makeidx}
\usepackage{tikz-network}
\usepackage[left=2cm,right=2cm,top=2cm,bottom=2cm]{geometry}
\usepackage{caption}
\usepackage{subcaption}
\usepackage[section]{placeins}
\usepackage{enumitem}
\usepackage{tikz}
\usepackage{stmaryrd}
\usepackage{upquote}
\usepackage{authblk}

\usepackage{multicol}

\usepackage{hyperref}
\hypersetup{colorlinks=true, citecolor={blue}, linkcolor=blue, filecolor=magenta, urlcolor=cyan}
\usetikzlibrary{arrows.meta, positioning, decorations.pathreplacing}
\usepackage{cleveref}
  
\newtheorem{theorem}{Theorem}[section]

\newtheorem{lemma}[theorem]{Lemma}

\newtheorem{proposition}[theorem]{Proposition}

\theoremstyle{definition}
\newtheorem{defi}[theorem]{Definition}

\numberwithin{equation}{section}

\newcommand{\Nl}{\mathbb{N}}
\newcommand{\Zl}{\mathbb{Z}}
\newcommand{\Rl}{\mathbb{R}}
\newcommand{\Ql}{\mathbb{Q}}

\newcommand{\Cl}{\mathbb{C}}
\newcommand{\eu}{\mathbf{e}_u}
\newcommand{\ev}{\mathbf{e}_v}
\newcommand{\eut}{\mathbf{e}_u^t}
\newcommand{\evt}{\mathbf{e}_v^t}

\newcommand{\ee}{\boldsymbol e}

\newcommand{\dist}{\operatorname{dist}}
\newcommand{\halfZ}{\Zl[1/2]}
\usepackage{titling}
\title{{Perfect state transfer in Grover walks on generalized Bethe trees}
}
\author{Koushik Bhakta\thanks{Corresponding author.\\ E-mail addresses: b.koushik@iitg.ac.in (K. Bhakta), b.bikash@iitg.ac.in (B. Bhattacharjya).} \  and Bikash Bhattacharjya\\
	Department of Mathematics\\
	Indian Institute of Technology Guwahati, India}

\date{}
\begin{document}
	\maketitle
	
	\vspace{-1cm}
	
	\begin{center}{\textbf{Abstract}}\end{center}

	\noindent
	We completely classify the generalized Bethe trees that exhibit perfect state transfer between vertex states in Grover walks. A generalized Bethe tree $B(b_0,\ldots,b_{n-1})$ is a rooted tree in which every vertex at level $i<n$ has $b_i$ children. We prove that such a tree exhibits perfect state transfer if and only if it is a path, $S_k(B(2,3))$ or $S_k(B(1,2,3))$, where $S_k(G)$ replaces every edge of $G$ by a path of length $k$. We also determine all transferring pairs and the minimum transfer times.  
	
	\vspace*{.5cm}
	
	\noindent 
	\textbf{Keywords.} Generalized Bethe trees; Grover walk; perfect state transfer; Chebyshev polynomial; Jacobi matrix.
	
	\noindent
	\textbf{Mathematics Subject Classifications:} 05C05; 05C50; 05C81; 81Q99.
	
	\section{Introduction and statement of results} 
	Quantum walks on graphs are quantum analogues of classical random walks~\cite{quantum_algorithm}. Their evolution is governed by unitary operators, and interference between different routes can lead to behaviour that does not occur in classical random walks. Quantum walks provide a framework for quantum algorithms~\cite{ambains} and for transferring quantum information through networks~\cite{christandl}. Quantum walks are broadly divided into continuous-time~\cite{state} and discrete-time~\cite{godsil_dqw} models according to their time evolution. In this paper, we consider the Grover walk, a widely studied discrete-time quantum walk that uses the Grover diffusion operator~\cite{hig3,first_grover}. In a Grover walk, the state assigns a complex amplitude to each arc of the graph. At each step, we apply the Grover diffusion operator at every vertex and then exchange the amplitudes of opposite arcs. Thus the  underlying graph determines the evolution operator~\cite{kubota_entangle}. We give the formal definition in Section~\ref{preli_grover}.

	Following Kubota and Segawa~\cite{kubota1}, we take the \emph{vertex state} at $u$ as the uniform superposition over all arcs with terminus $u$. We say that \emph{perfect state transfer} occurs from $u$ to $v$ if a Grover walk starting in the vertex state at $u$ reaches the vertex state at $v$ with probability $1$ after some time. This notion is motivated by the problem of transmitting quantum information through the natural evolution of a spin network~\cite{christandl}.  We ask how the structure and spectrum of a graph determine whether it exhibits such transfer between distinct vertices.

	Perfect state transfer in Grover walks is rare, and the known examples are very symmetric.  Barr et al.~\cite{barr_pst} studied perfect state transfer on variants of cycles. Zhan~\cite{zhan1} constructed an infinite family of circulant graphs of valency $4$. Using Chebyshev polynomials, Kubota and Segawa~\cite{kubota1} obtained a necessary spectral condition for perfect state transfer and characterized its occurrence on complete multipartite graphs with equal part sizes. Kubota and Yoshino~\cite{kubota2} later classified circulant graphs of valency at most $4$ that exhibit perfect state transfer.  For finite abelian groups $\Gamma$, Bhakta and Bhattacharjya~\cite{bhakta_abelian} determined the pairs $(\Gamma,S)$ for which the Cayley graph $\operatorname{Cay}(\Gamma,S)$  exhibits perfect state transfer. Bhakta, Bhattacharjya, and Cao~\cite{bhakta_dihedral} obtained the corresponding classification for dihedral groups.

	Most of the families discussed above consist of vertex-transitive graphs. On such a graph, perfect state transfer between two vertex states implies that the Grover walk is periodic~\cite{bhakta3}. So a natural question is what happens when this symmetry is gone. Trees are the opposite case. A tree is the cheapest connected network, with the fewest possible edges, so it is a natural model for a network of qubits. For continuous-time walks the answer on trees is discouraging. For the walk defined by the Laplacian matrix, Coutinho and Liu~\cite{coutinho_laplacian_tree} showed that only the paths on two vertices exhibits perfect state transfer and conjectured that paths on two and three vertices are the only such trees in the adjacency model. Coutinho, Juliano, and Spier~\cite{coutinho_tree} later proved this conjecture.
	
	The Grover walk behaves differently. At a vertex of degree two, the Grover coin simply swaps the two directions~\cite{mixed2}. So the amplitude runs straight along a path, and paths of every length exhibit perfect state transfer. The continuous-time barrier is gone, at least for paths. This leads to a clear question: which trees with branching can still transfer a state perfectly?
	
	General trees are hard to control, so we start with a class that has branching but is still easy to analyze. A \emph{generalized Bethe tree} is a rooted tree in which all vertices at the same distance from the root have the same degree. The degree can change from one level to the next. Paths and stars are examples, and every such tree with more than two vertices is nonregular. The symmetry inside each level allows an explicit spectral analysis. Kubota, Segawa, Taniguchi, and Yoshie~\cite{bethetree} classified the generalized Bethe trees whose Grover walks are periodic. Their list contains the paths and some families of subdivided trees. This is a good start, but it does not answer our question, for two reasons. Periodicity says that every state returns to itself. It does not say that a state can reach a different vertex, and perfect state transfer only needs the transferring state to be periodic, not the whole walk. So the periodic trees cannot simply be read as the answer.
	
	In this paper we study the problem directly, without assuming that the whole walk is periodic. To the best of our knowledge, this is the first study of perfect state transfer between vertex states in Grover walks on trees. We give a complete answer for generalized Bethe trees. Besides the paths, perfect state transfer occurs on two infinite families of nonregular trees with branching. To state the result, let $B(b_0,\ldots,b_{n-1})$ be the generalized Bethe tree in which every vertex at level $i<n$ has $b_i$ children. For a positive integer $k$, let $S_k(G)$ be the graph obtained by replacing every edge of $G$ by a path of length $k$. Our main result is the following.

	\begin{theorem}\label{pst_gbt}
		A generalized Bethe tree with at least two vertices exhibits
		perfect state transfer in Grover walks if and only if it is
		isomorphic to $P_{n+1}$,
		$S_k(B(2,3))$, or $S_k(B(1,2,3))$ for positive integers $n$ and $k$.
	\end{theorem}
		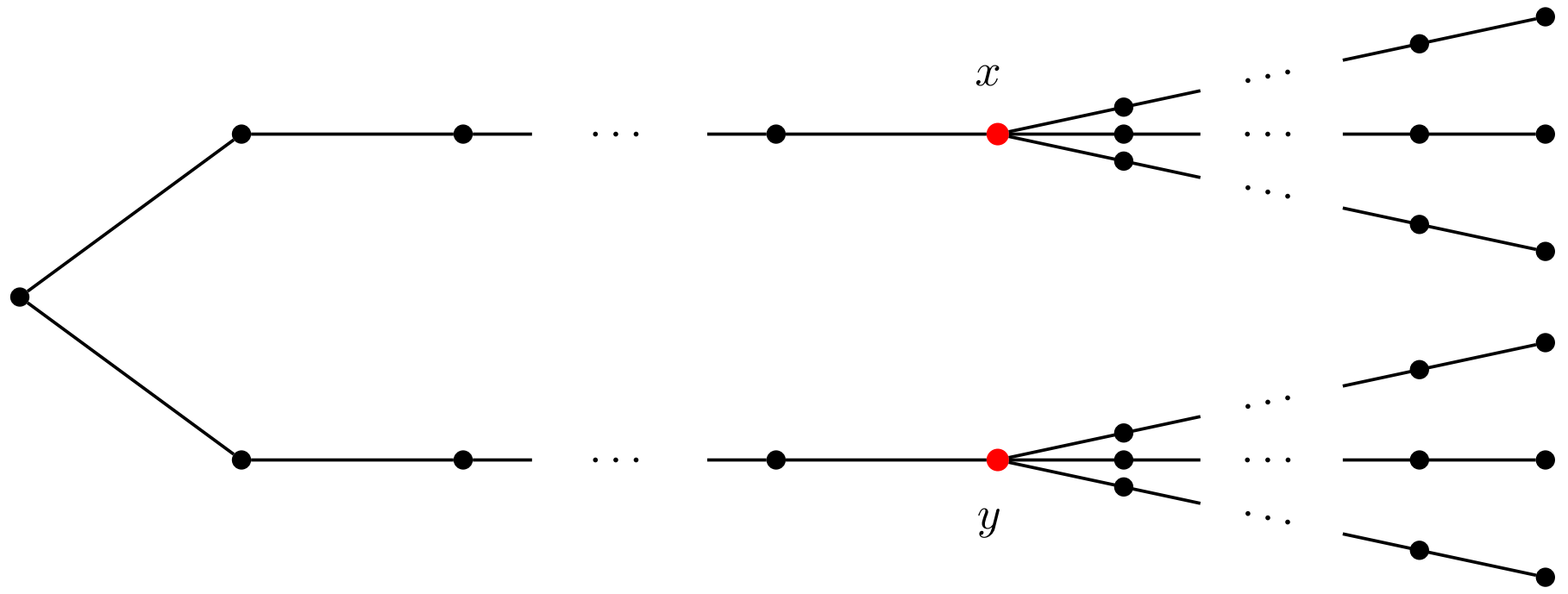
\begin{figure}[h!]
		\centering
		\begin{tikzpicture}[
			x=1cm,y=1cm,font=\normalsize,
			v/.style={circle,fill=black,draw=black,inner sep=0pt,minimum size=3.8pt},
			p/.style={circle,fill=red,draw=red,line width=.9pt,inner sep=0pt,minimum size=4pt},
			e/.style={line width=.7pt},
			continuation/.style={midway,sloped,anchor=base,inner sep=0pt,outer sep=0pt},
			every label/.style={inner sep=2pt,font=\normalsize}
			]
			
			\node[v] (r) at (0,0) {};
			\node[v] (x1) at (1.7,1.25) {};
			\node[v] (x2) at (3.4,1.25) {};
			\node[v] (xm) at (5.8,1.25) {};
			\node[p] (x) at (7.5,1.25) {};
			\node[anchor=south east] at (7.65,1.5) {$x$};
			
			\node[v] (y1) at (1.7,-1.25) {};
			\node[v] (y2) at (3.4,-1.25) {};
			\node[v] (ym) at (5.8,-1.25) {};
			\node[p] (y) at (7.5,-1.25) {};
			\node[anchor=north east] at (7.65,-1.5) {$y$};
			
			\draw[e] (r)--(x1)--(x2);
			\draw[e] (r)--(y1)--(y2);
			\foreach \a/\b in {x2/xm,y2/ym}{
				\draw[e] (\a)--($(\a)!.22!(\b)$);
				\path ($(\a)!.22!(\b)$)--node[continuation]
				{$\raisebox{-\fontdimen22\textfont2}{$\cdots$}$}
				($(\a)!.78!(\b)$);
				\draw[e] ($(\a)!.78!(\b)$)--(\b);
			}
			\draw[e] (xm)--(x);
			\draw[e] (ym)--(y);

			\foreach \s/\yy in {1/2.15,2/1.25,3/.35}{
				\node[v] (a\s) at (11.7,\yy) {};
				\coordinate (aa\s) at ($(x)!.23!(a\s)$);
				\coordinate (ab\s) at ($(x)!.77!(a\s)$);
				\node[v] at (aa\s) {};
				\node[v] at (ab\s) {};
				\draw[e] (x)--($(x)!.37!(a\s)$);
				\path ($(x)!.37!(a\s)$)--node[continuation]
				{$\raisebox{-\fontdimen22\textfont2}{$\cdots$}$}
				($(x)!.63!(a\s)$);
				\draw[e] ($(x)!.63!(a\s)$)--(a\s);
			}
			
			\foreach \s/\yy in {1/-.35,2/-1.25,3/-2.15}{
				\node[v] (b\s) at (11.7,\yy) {};
				\coordinate (ba\s) at ($(y)!.23!(b\s)$);
				\coordinate (bb\s) at ($(y)!.77!(b\s)$);
				\node[v] at (ba\s) {};
				\node[v] at (bb\s) {};
				\draw[e] (y)--($(y)!.37!(b\s)$);
				\path ($(y)!.37!(b\s)$)--node[continuation]
				{$\raisebox{-\fontdimen22\textfont2}{$\cdots$}$}
				($(y)!.63!(b\s)$);
				\draw[e] ($(y)!.63!(b\s)$)--(b\s);
			}
		\end{tikzpicture}
		\vspace{2mm}
		
		\caption{$S_k(B(2,3))$, with the unique transferring pair $\{x,y\}$.}
		\label{fig:pst-b23}
	\end{figure}
		
	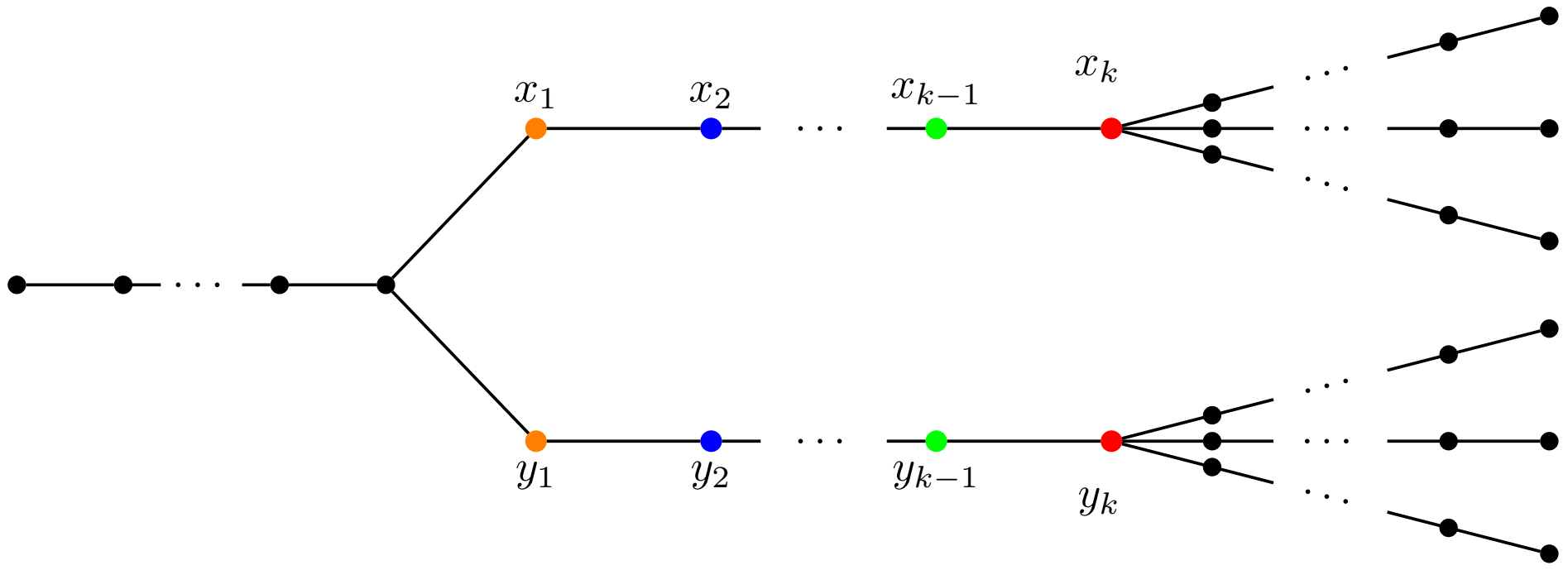
\begin{figure}[h!]
		\centering
		\begin{tikzpicture}[
			x=1cm,y=1cm,font=\normalsize,
			v/.style={circle,fill=black,draw=black,inner sep=0pt,minimum size=3.8pt},
			p/.style={circle,fill=orange,draw=orange,line width=.9pt,inner sep=0pt,minimum size=4pt},
			q/.style={circle,fill=blue,draw=blue,line width=.9pt,inner sep=0pt,minimum size=4pt},
			r/.style={circle,fill=green,draw=green,line width=.9pt,inner sep=0pt,minimum size=4pt},
			s/.style={circle,fill=red,draw=red,line width=.9pt,inner sep=0pt,minimum size=4pt},
			e/.style={line width=.7pt},
			continuation/.style={midway,sloped,anchor=base,inner sep=0pt,outer sep=0pt},
			every label/.style={inner sep=2pt,font=\normalsize}
			]
			
			\node[v] (z1) at (.85,0) {};
			\node[v] (z2) at (1.7,0) {};
			\node[v] (zm) at (2.95,0) {};
			\node[v] (c) at (3.8,0) {};
			\draw[e] (z1)--(z2)--(2,0);
			\path (2,0)--node[continuation]
			{$\raisebox{-\fontdimen22\textfont2}{$\cdots$}$} (2.65,0);
			\draw[e] (2.65,0)--(zm)--(c);

			\node[p,label=above:$x_1$] (x1) at (5,1.25) {};
			\node[q,label=above:$x_2$] (x2) at (6.4,1.25) {};
			\node[r,label=above:$x_{k-1}$] (xm) at (8.2,1.25) {};
			\node[s] (x) at (9.6,1.25) {};
			\node[anchor=south east] at (9.8,1.5) {$x_k$};
			
			\node[p,label=below:$y_1$] (y1) at (5,-1.25) {};
			\node[q,label=below:$y_2$] (y2) at (6.4,-1.25) {};
			\node[r,label=below:$y_{k-1}$] (ym) at (8.2,-1.25) {};
			\node[s] (y) at (9.6,-1.25) {};
			\node[anchor=north east] at (9.8,-1.5) {$y_k$};
			
			\draw[e] (c)--(x1)--(x2);
			\draw[e] (c)--(y1)--(y2);
			\foreach \a/\b in {x2/xm,y2/ym}{
				\draw[e] (\a)--($(\a)!.22!(\b)$);
				\path ($(\a)!.22!(\b)$)--node[continuation]
				{$\raisebox{-\fontdimen22\textfont2}{$\cdots$}$}
				($(\a)!.78!(\b)$);
				\draw[e] ($(\a)!.78!(\b)$)--(\b);
			}
			\draw[e] (xm)--(x);
			\draw[e] (ym)--(y);
			
			\foreach \s/\yy in {1/2.15,2/1.25,3/.35}{
				\node[v] (a\s) at (13.1,\yy) {};
				\node[v] at ($(x)!.23!(a\s)$) {};
				\node[v] at ($(x)!.77!(a\s)$) {};
				\draw[e] (x)--($(x)!.37!(a\s)$);
				\path ($(x)!.37!(a\s)$)--node[continuation]
				{$\raisebox{-\fontdimen22\textfont2}{$\cdots$}$}
				($(x)!.63!(a\s)$);
				\draw[e] ($(x)!.63!(a\s)$)--(a\s);
			}
			\foreach \s/\yy in {1/-.35,2/-1.25,3/-2.15}{
				\node[v] (b\s) at (13.1,\yy) {};
				\node[v] at ($(y)!.23!(b\s)$) {};
				\node[v] at ($(y)!.77!(b\s)$) {};
				\draw[e] (y)--($(y)!.37!(b\s)$);
				\path ($(y)!.37!(b\s)$)--node[continuation]
				{$\raisebox{-\fontdimen22\textfont2}{$\cdots$}$}
				($(y)!.63!(b\s)$);
				\draw[e] ($(y)!.63!(b\s)$)--(b\s);
			}
		\end{tikzpicture}
		\vspace{2mm}
		
		\caption{$S_k(B(1,2,3))$, with transferring pairs $\{x_i,y_i\}$ for $1\leq i\leq k$.}
		\label{fig:pst-b123}
	\end{figure}
	We now describe all transferring pairs and their minimum transfer times. The two nonpath families are illustrated in Figures~\ref{fig:pst-b23} and~\ref{fig:pst-b123}. Label the vertices of $P_{n+1}$ as $v_0,\ldots,v_n$ in path order. In each nonpath family, let $x,y$ be the two vertices of degree $4$. For $S_k(B(1,2,3))$, let $c$ be the unique vertex of degree $3$, and write the paths from $c$ to $x$ and $y$ as $c=x_0,x_1,\ldots,x_k=x$ and  $c=y_0,y_1,\ldots,y_k=y$.
	
	\begin{proposition}\label{prop:pst-pairs}
		For each graph in Theorem~\ref{pst_gbt}, each pairs of distinct vertices exhibiting perfect state transfer and their minimum transfer times are given in Table~\ref{tab:pst-pairs}.
	\end{proposition}
	\begin{table}[ht]
		\centering
		\small
		\setlength{\tabcolsep}{25pt}
		\renewcommand{\arraystretch}{1.5}
		\begin{tabular}{@{}ccc@{}}
			\hline
			Graph & Minimum time & Transferring pairs \\
			\hline
			$P_{n+1}$ 
			& $n$
			& $\{v_j,v_{n-j}\}$  $(0\leq j\leq n-1-j)$ \\[2pt]
			
			$S_k(B(2,3))$  
			& $6k$
			& $\{x,y\}$ \\[2pt]
			
			$S_k(B(1,2,3))$  
			& $6k$
			& $\{x_j,y_j\}$  $(1\leq j\leq k)$ \\
			\hline
		\end{tabular}
		\caption{Minimum transfer times and all transferring pairs.}
		\label{tab:pst-pairs}
	\end{table}

	The graph $S_k(B(2,3))$ has exactly one transferring pair, whereas $S_k(B(1,2,3))$ has exactly $k$ transferring pairs. In the second family, the pairs consist of vertices at equal distances from $c$ along the two paths to $x$ and $y$.

	\section{Preliminaries}
	Throughout the paper, $G$ denotes a simple connected graph with vertex set $V(G)$ and edge set $E(G)$. Each edge $\{u,v\}$ is an unordered pair of vertices and determines two arcs, $(u,v)$ and $(v,u)$. The symmetric arc set of $G$ is $\mathcal{A}(G)=\{(u,v), (v,u):\{u,v\}\in E(G)\}$. For an arc $a=(u,v)$, we call $u$ the \emph{origin} and $v$ the \emph{terminus} of $a$, and write $o(a)=u$ and $t(a)=v$. The \emph{inverse arc} of $a$ is $a^{-1}=(v,u)$.
	\subsection{Grover walks and perfect state transfer}\label{preli_grover}
	We now define the matrices associated with the Grover walk on a graph $G$. We refer to \cite{kubota1} for more information about the matrices. The \emph{boundary matrix} $N:=N(G)\in\mathbb{C}^{V(G)\times\mathcal{A}(G)}$ of $G$ has entries $N_{u,a}=\delta_{u,t(a)}/\sqrt{\deg(u)}$, where $\delta_{x,y}$ denote the Kronecker delta. It can be checked that $NN^*=I$. The \emph{shift matrix} $R:=R(G)\in\mathbb{C}^{\mathcal{A}(G)\times\mathcal{A}(G)}$ has entries $R_{a,b}=\delta_{a,b^{-1}}$. The \emph{time evolution matrix} $U:=U(G)\in\mathbb{C}^{\mathcal{A}(G)\times\mathcal{A}(G)}$, also called the \emph{Grover transition matrix}, is defined by
	\[
	U=R(2N^*N-I),
	\]
	where $N^*$ denotes the conjugate transpose of $N$ and $I$ is the identity matrix on $\mathbb{C}^{\mathcal{A}(G)}$. The discrete-time quantum walk determined by $U$ is called the \emph{Grover walk} on the graph  $G$. The Grover walks are a special case of bipartite walks~\cite{chen}.
	
	A state of the Grover walk on $G$ is a unit vector in $\Cl^{\mathcal{A}(G)}$. We say that \emph{perfect state transfer} occurs from a state $\Phi_1$ to a distinct state $\Phi_2$ at time $\tau\in\Nl$ if $U^\tau\Phi_1=\gamma\Phi_2$ for some $\gamma\in\Cl$ with $|\gamma|=1$.

	We focus on the transfer of states localized at individual vertices of a graph. We associate with each vertex $u$ the \emph{vertex-type state} $\Phi_u:=N^*\eu$, where $\eu$ is the standard basis vector at $u$. For the motivation to restrict states under consideration to vertex-type states, we refer the reader to \cite{kubota1}.  Throughout this paper, we focus on  vertex-type states.  
	\begin{defi}
		A graph exhibits \emph{perfect state transfer} from a vertex $u$ to a vertex $v$ at time $\tau$ if there exists a complex number $\gamma$ with $|\gamma|=1$ such that $U^\tau \Phi_u=\gamma \Phi_v$.
	\end{defi}
	Perfect state transfer is always understood between distinct vertices. We say that a graph $G$ exhibits \emph{perfect state transfer} if there exist vertices $u\neq v$ in $G$ and a positive integer $\tau$ such that perfect state transfer occurs from $u$ to $v$ at time $\tau$.

	Perfect state transfer in Grover walks on graphs can be studied effectively using the discriminant matrix and Chebyshev polynomials (see \cite{kubota2}). The \emph{discriminant} $P:=P(G)\in\mathbb{C}^{V(G)\times V(G)}$ of $G$ is defined by $P=NRN^*$. It is similar to the transition matrix of the isotropic random walk on $G$. Indeed, $P$ is the normalized adjacency matrix of $G$, that is, $P=D^{-1/2}AD^{-1/2}$, where $A$ is the adjacency matrix and $D$ is the diagonal degree matrix of $G$ \cite{mixed1}.
	
	The \emph{Chebyshev polynomials of the first kind}, denoted by $T_n(x)$, are defined by the recurrence relation $T_0(x)=1$, $T_1(x)=x$, and $T_n(x)=2xT_{n-1}(x)-T_{n-2}(x)$ for $n\geq2$. It is well known that $T_n(\cos\theta)=\cos(n\theta)$. Consequently, for $n\geq1$ and $x\in[-1,1]$, we have $|T_n(x)|=1$ if and only if $x=\cos(j\pi/n)$ for some integer $j$.
	\begin{lemma}[{\cite[Lemma~3.1]{kubota1}}]\label{P_and_U}
		Let $G$ be a finite simple graph with the transition matrix $U$ and the discriminant matrix $P$. Then $NU^kN^*=T_k(P)$ for any nonnegative integer $k$.
	\end{lemma}

	Let $\sigma(P)$ denote the set of all distinct eigenvalues of $P$. Since $P$ is symmetric, it has spectral decomposition. Let the spectral decompostion of $P$ be 
	\begin{equation}\label{crb1749}
	P=\sum_{\mu\in\sigma(P)}\mu E_\mu,
	\end{equation}
	where $E_\mu$ is the orthogonal projection onto the $\mu$-eigenspace. The eigenvalue support of a vertex $u$ is $\Theta_P(u)=\{\mu\in\sigma(P):E_\mu\eu\neq0\}$. We say two vertices $u$ and $v$ are strongly cospectral if $E_\mu\eu=\pm E_\mu\ev$ for each $\mu\in\sigma(P)$.

	\begin{lemma}[{\cite[Theorem~6.5]{kubota2}}]\label{pre_lemma1}
	A graph exhibits perfect state transfer from $u$ to $v$ at time $\tau$ if and only if $T_\tau(P)\eu=\ev$.
    \end{lemma}
 \begin{proof}
 	Suppose perfect state transfer occurs from $u$ to $v$ at time $\tau$, that is, $U^\tau\Phi_u=\gamma\Phi_v$ for some $|\gamma|=1$. Multiplying $N$ and using Lemma~\ref{P_and_U}, we obtain $T_\tau(P)\eu=\gamma\ev$. Let $q\in\Cl^{V(G)}$ with $q_w=\sqrt{\deg(w)}$ for each $w\in V(G)$. Since $Pq=q$ and $T_\tau(1)=1$, we have $T_\tau(P)q=q$. Thus, by the symmetry of $T_\tau(P)$,
 	\begin{equation}\label{crb1508}
 	\sqrt{\deg(u)}=q^t\eu=q^tT_\tau(P)\eu
 	=\gamma\sqrt{\deg(v)}.
 	\end{equation}
 	Since $|\gamma|=1$, we must have $\gamma=1$. Hence $T_\tau(P)\eu=\ev$.
 	
 	Conversely, suppose that $T_\tau(P)\eu=\ev$. By Lemma~\ref{P_and_U}, $\langle\Phi_v,U^\tau\Phi_u\rangle=\langle\ev,T_\tau(P)\eu\rangle=1$. Since $\Phi_v$ and $U^\tau\Phi_u$ are unit vectors, $U^\tau\Phi_u=\Phi_v$. Thus perfect state transfer occurs from $u$ to $v$ at time $\tau$.
 \end{proof}
 
	\begin{lemma}[\cite{bhakta3,kubota2}]\label{pre_lemma2}
	Let $G$ be a graph exhibits perfect state transfer from $u$ to $v$ at time $\tau$. Then the following hold.
	\begin{enumerate}
		\item[(i)] $T_\tau(P)\eu=\ev$ and $T_\tau(P)\ev=\eu$.
		\item[(ii)] The vertices $u$ and $v$ are strongly cospectral with respect to $P$.
		\item[(iii)] $\deg(u)=\deg(v)$.
		\item[(iv)] Every automorphism of $G$ that fixes $u$ also fixes $v$.
		\item[(v)] For each $\mu\in\Theta_P(u)$,  $T_\tau(\mu)=\pm1$. In particular, $2\mu$ is an algebraic integer.
	\end{enumerate}
\end{lemma}
 \begin{proof}
 	By Lemma~\ref{pre_lemma1}, $T_\tau(P)\eu=\ev$. Since $NN^*=I$ and $U$ is unitary, for each $x\in V(G)$ we have $\|T_\tau(P)\ee_x\|=\|NU^\tau N^*\ee_x\|\leq \|U^\tau N^*\ee_x\|=\|N^*\ee_x\|=1$.	By the symmetry of $T_\tau(P)$, we have $\eut T_\tau(P)\ev=\evt T_\tau(P)\eu=1$. Thus  $T_\tau(P)\ev=\eu$. This proves \textup{(i)}.

 	By \eqref{crb1749} and \textup{(i)}, we find  $E_\mu\ev=T_\tau(\mu)E_\mu\eu$ and $E_\mu\eu=T_\tau(\mu)E_\mu\ev$ for each $\mu\in\sigma(P)$.  If $\mu\in\Theta_P(u)$, then $E_\mu\eu\neq0$, and so $T_\tau(\mu)^2=1$. Hence $E_\mu\ev=\pm E_\mu\eu$. If $\mu\notin\Theta_P(u)$, $E_\mu\ev=E_\mu\eu=0$. Thus $u$ and $v$ are strongly cospectral, proving \textup{(ii)}.
 	
 	The degree equality in \textup{(iii)} follows from \eqref{crb1508}.
 	
 	For \textup{(iv)}, let $\sigma$ be an automorphism of $G$ that fixes $u$, and let $Q$ be its permutation matrix. Then $Q$ commutes with $P$ and hence with $T_\tau(P)$. Therefore $\ee_{\sigma(v)}=Q\ev=QM\eu=MQ\eu=M\eu=\ev$. 	Hence $\sigma(v)=v$.
 	
 	Finally, let $\mu\in\Theta_P(u)$. We have shown that $\alpha:=T_\tau(\mu)\in\{-1,1\}$. By the Chebyshev recurrence, $2T_\tau(x/2)-2\alpha$ is a monic polynomial in $\Zl[x]$. Since $2\mu$ is a root of this polynomial, it is an algebraic integer. This proves \textup{(v)}.
 \end{proof}
 
	\subsection{Generalized Bethe trees}
	
	A \emph{rooted tree} is a tree with a chosen vertex, called  the \emph{root}. The \emph{level} of a vertex is its distance from the root. Note that the root lies at level $0$. Every vertex other than the root has a unique neighbour in the preceding level, called its \emph{parent}, and the neighbours of a vertex in the next level  are called its \emph{children}.
	
	\begin{defi}
		For positive integers $b_0,\ldots,b_{n-1}$, the \emph{generalized Bethe tree} $B(b_0,\ldots,b_{n-1})$ is the rooted
		tree in which every vertex at level $i$ has $b_i$ children for $0\leq i<n$, and the vertices at level $n$ have no children.
	\end{defi}

\begin{figure}[h!]
		\centering
		\begin{tikzpicture}[x=1cm,y=1cm,font=\normalsize,v/.style={circle,fill=black,inner sep=0pt,minimum size=4pt},e/.style={line width=.7pt}]
			
			\foreach \i/\N in {0/1,1/2,2/6,3/6}{
				\fill[gray!10,rounded corners=3pt]
				({3.4*\i-.4},-2.9) rectangle ({3.4*\i+.4},2.9);
				\node at ({3.4*\i},3.3) {$C_{\i}$};
				\node at ({3.4*\i},-3.3) {$c_{\i}=\N$};
			}
			
			\node[v] (r) at (0,0) {};
			\node[v] (a) at (3.4,1.5) {};
			\node[v] (b) at (3.4,-1.5) {};
			\draw[e] (r)--(a) (r)--(b);
			
			\foreach \j/\yy in {1/2.5,2/1.5,3/.5,4/-.5,5/-1.5,6/-2.5}{
				\node[v] (u\j) at (6.8,\yy) {};
				\node[v] (w\j) at (10.2,\yy) {};
				\draw[e] (u\j)--(w\j);
			}
			\foreach \j in {1,2,3}{
				\draw[e] (a)--(u\j);
			}
			\foreach \j in {4,5,6}{
				\draw[e] (b)--(u\j);
			}
			
			\node at (1.7,0) {$b_0=2$};
			\node at (5.1,0) {$b_1=3$};
			\node at (8.5,0) {$b_2=1$};
		\end{tikzpicture}
		
		\caption{The generalized Bethe tree $B(2,3,1)$.}
		\label{fig:bethe-levels}
	\end{figure}
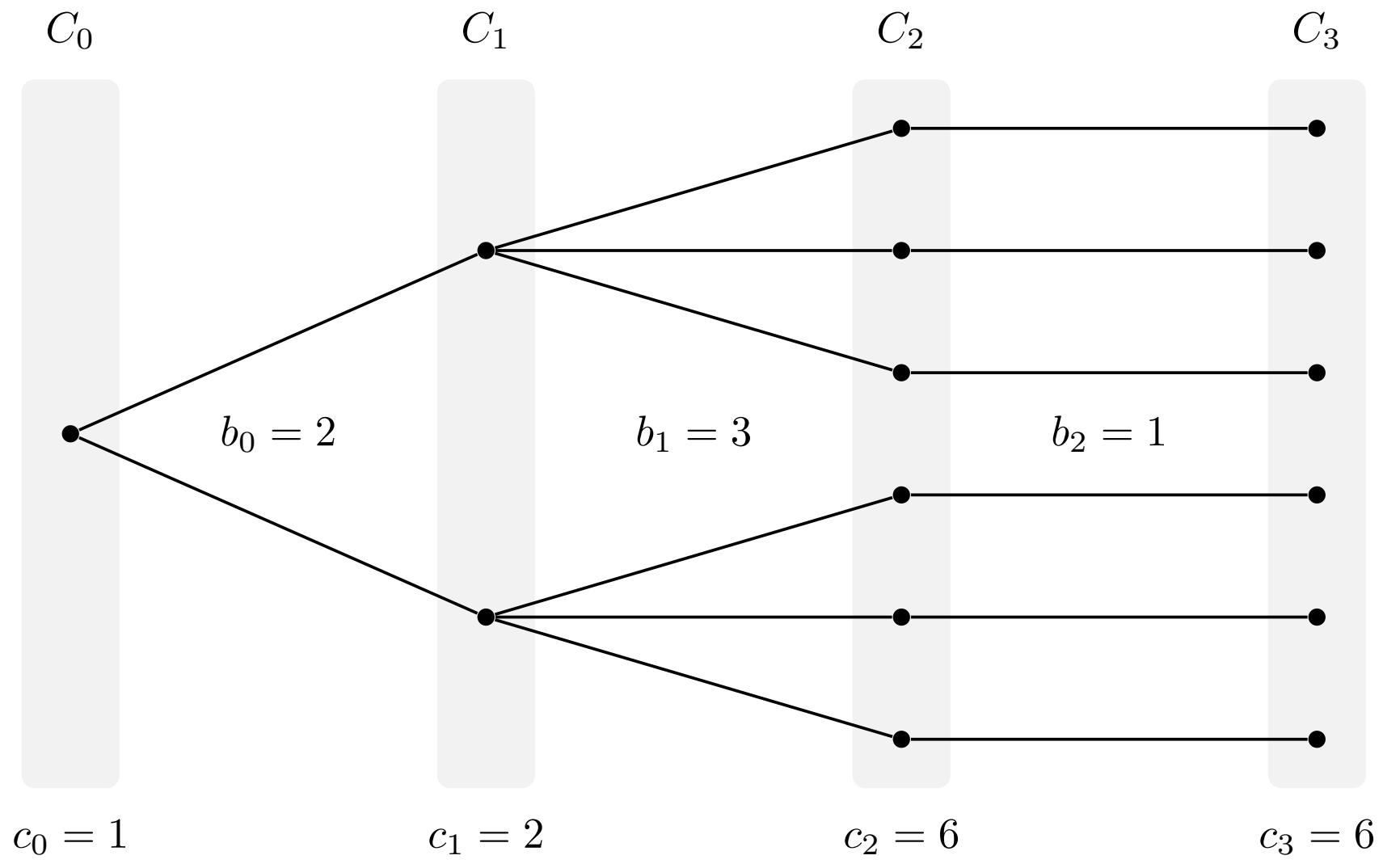

	Throughout the paper, $B$ denotes the generalized Bethe tree $B(b_0,\ldots,b_{n-1})$, unless stated otherwise.  For $0\leq i<n$, the positive integer $b_i$ is called the \emph{branching number} at level $i$. We call level $i$ a \emph{branching level} if $b_i\geq2$. For $0\leq i\leq n$, let $C_i$ denote the set of vertices at level $i$, and put $c_i=|C_i|$. Then $c_0=1$ and $c_{i+1}=c_i b_i$ for $0\leq i<n$. Let $d_i$ be the common degree of the vertices in level $i$. Then $d_0=b_0$, $d_i=b_i+1$ for  $1\leq i<n$, and $d_n=1$. Figure~\ref{fig:bethe-levels} shows the generalized Bethe tree $B(2,3,1)$ with its levels $C_0$, $C_1$, $C_2$ and $C_3$.

	Let $V$ be the vertex set of $B(b_0,\ldots,b_{n-1})$. For $0\leq i\leq n$, define $\Psi_i=(1/\sqrt{c_i})\mathbf{1}_i\in\Cl^{|V|}$, where $\mathbf{1}_i$ is the indicator vector of $C_i$, that is, $\mathbf{1}_i(v)=1$ if $v\in C_i$ and $0$ otherwise. The vectors $\Psi_0,\ldots,\Psi_n$ form an orthonormal set in
	$\Cl^{|V|}$. We denote their span by $\mathcal{R}=\operatorname{span}\{\Psi_0,\ldots,\Psi_n\}$.
	
	For $0\leq i\leq n-1$, put $\beta_i=\sqrt{b_i/(d_i d_{i+1})}$. Using $c_{i+1}=c_i b_i$, a direct calculation gives $P\Psi_0=\beta_0\Psi_1$, $P\Psi_i=\beta_{i-1}\Psi_{i-1}+\beta_i\Psi_{i+1}$ for $1\leq i\leq n-1$, and $P\Psi_n=\beta_{n-1}\Psi_{n-1}$. These calculations can be found in \cite[Lemma~2.4]{bethetree}.
	
	Hence $\mathcal R$ is invariant under $P$. With respect to the ordered orthonormal basis $\{\Psi_0,\ldots,\Psi_n\}$, the restriction $P|_{\mathcal R}$ is represented by the symmetric tridiagonal matrix
	\begin{equation}\label{eq:radial-matrix}
		J_0:=
		\begin{pmatrix}
			0 & \beta_0 & & & \\
			\beta_0 & 0 & \beta_1 & & \\
			& \beta_1 & 0 & \ddots & \\
			& & \ddots & \ddots & \beta_{n-1} \\
			& & & \beta_{n-1} & 0
		\end{pmatrix},
	\end{equation}
	where all blank entries are zero. 
	
	For $1\leq i\leq n$, let $J_i$ be the matrix obtained from $J_0$ by deleting its first $i$ rows and first $i$ columns. The rows and columns of $J_i$ are indexed by the levels $i,i+1,\ldots,n$. The matrix $J_i$ has order $(n-i+1)\times(n-i+1)$.
	
	\begin{lemma}
		For $1\leq i\leq n$,
		\begin{equation}\label{det_J_i}
			\det(I-J_i)=\frac{1}{\prod_{\ell=i}^{n-1}(b_\ell+1)},
		\end{equation}
		where an empty product is $1$.
	\end{lemma}
	
	\begin{proof}
		Put $\Delta_i=\det(I-J_i)$. We prove the identity by induction on $n-i$. For $i=n$, we have $J_n=(0)$ and $\Delta_n=1$. If $n\geq2$, then $\Delta_{n-1}=1-\beta_{n-1}^2={1}/({b_{n-1}+1})$. Now let $1\leq i\leq n-2$ and assume that
		\eqref{det_J_i} holds for $i+1$ and $i+2$. Expanding $\det(I-J_i)$ along the first row gives $\Delta_i=\Delta_{i+1}-\beta_i^2\Delta_{i+2}$. By the induction hypothesis,
		\begin{align*}
			\Delta_i
			&=\frac{1}{\prod_{\ell=i+1}^{n-1}(b_\ell+1)}
			-\frac{b_i}{(b_i+1)(b_{i+1}+1)}
			\frac{1}{\prod_{\ell=i+2}^{n-1}(b_\ell+1)}\\
			&=\frac{1}{\prod_{\ell=i+1}^{n-1}(b_\ell+1)}
			\left(1-\frac{b_i}{b_i+1}\right)\\
			&=\frac{1}{\prod_{\ell=i}^{n-1}(b_\ell+1)}.
		\end{align*}
		This proves \eqref{det_J_i}.
	\end{proof}

	For a vertex $w\in C_k$ and $k\leq s\leq n$, let $D_s(w)$ be the set of vertices in $C_s$ whose path from the root passes through $w$. In particular, $D_k(w)=\{w\}$. The set $D_s(w)$ is also known as the descendants of $w$ at level $s$ (see \cite[Section~2]{bethetree}). If $w_1,\ldots,w_{b_j}$ are the children of a vertex $w\in C_j$, then each set $D_s(w_\ell)$ has the vertices
	\[m_{j,s}:=\prod_{t=j+1}^{s-1}b_t\quad\text{for}\quad j+1\leq s\leq n,\]
	with the convention that an empty product is $1$. In particular, $m_{j,j+1}=1$ and $m_{j,s+1}=b_s m_{j,s}$ for $j+1\leq s<n$.

	The vectors $\Psi_0,\ldots,\Psi_n$ are constant on their respective levels. We now construct vectors whose entries sum to zero on a level. For a subset $X\subseteq V$, let $\mathbf{1}_X$ denote its  indicator vector.
	
	Let $0\leq j<n$ with $b_j\geq2$. For each $w\in C_j$, denote its children by $w_1,\ldots,w_{b_j}$. Let $a=(a_1,\ldots,a_{b_j})$ be a unit vector in $\Rl^{b_j}$ with $\sum_{\ell=1}^{b_j}a_\ell=0$. For $j+1\leq s\leq n$, define
	\[
	\Phi_s(w,a):=\frac{1}{\sqrt{m_{j,s}}}
	\sum_{\ell=1}^{b_j}a_\ell\mathbf{1}_{D_s(w_\ell)}.
	\]
	Thus $\Phi_s(w,a)$ takes the value $a_\ell/\sqrt{m_{j,s}}$ on the
	level-$s$ descendants of $w_\ell$ and is zero elsewhere. Define
	\[
	\mathcal W(w,a):= \operatorname{span}\{\Phi_{j+1}(w,a),\ldots,\Phi_n(w,a)\}.
	\]
	
	At each branching level $j$ and each $w\in C_j$, we use the following $b_j-1$ choices of $a$,
	\[
	a^{(q)}
	=\frac{1}{\sqrt{q(q+1)}}
	\bigl(\underbrace{1,\ldots,1}_{q},-q,0,\ldots,0\bigr),
	\qquad 1\leq q\leq b_j-1.
	\]
	Each $a^{(q)}$ is a unit vector with zero sum and $\langle a^{(p)},a^{(q)}\rangle=0$ for $p\neq q$. Thus each $w\in C_j$ gives $b_j-1$ subspaces $\mathcal W(w,a^{(1)}),\ldots,\mathcal W(w,a^{(b_j-1)})$.
	
    \begin{lemma}\label{invariant-orthogonal}
	With the preceding notations, $\Phi_{j+1}(w,a),\ldots,\Phi_n(w,a)$ is an orthonormal basis of the $P$-invariant subspace $\mathcal W(w,a)$, and $P|_{\mathcal W(w,a)}$ is represented by $J_{j+1}$ in this basis. Moreover, $\mathcal R$ and the subspaces $\mathcal W(w,a^{(q)})$, for all branching levels $j$, $w\in C_j$ and $1\leq q\leq b_j-1$, are mutually orthogonal.
    \end{lemma}
	\begin{proof}
		Fix $w\in C_j$ and $j+1\leq s\leq n$. The sets $D_s(w_1),\ldots,D_s(w_{b_j})$ are disjoint and each has $m_{j,s}$ vertices, so for any $a,a'$,
		\begin{equation}\label{eq:phi-inner}
			\langle\Phi_s(w,a),\Phi_s(w,a')\rangle=\sum_{\ell=1}^{b_j}a_\ell a'_\ell=\langle a,a'\rangle.
		\end{equation}
		Taking $a'=a$ gives $\|\Phi_s(w,a)\|=1$. Since $\Phi_s(w,a)(v)=0$ for $v\notin C_s$, we have $\langle\Phi_s(w,a),\Phi_t(w,a')\rangle=0$ for $s\neq t$. Hence $\Phi_{j+1}(w,a),\ldots,\Phi_n(w,a)$ is an orthonormal basis of $\mathcal W(w,a)$.
		
		Fix $a$ and write $\Phi_s=\Phi_s(w,a)$. We compute $(P\Phi_s)(v)$ in the following cases.

		\noindent(i) Let $s<n$ and $v\in D_{s+1}(w_\ell)$. The only neighbour of $v$ at level $s$ is its parent. Since $m_{j,s+1}=b_s m_{j,s}$,
		\begin{equation}\label{eq:P-down}
			(P\Phi_s)(v)
			=\frac{a_\ell}{\sqrt{m_{j,s}d_sd_{s+1}}}
			=\beta_s\Phi_{s+1}(v).
		\end{equation}
		
		\noindent(ii) Let $s>j+1$ and $v\in D_{s-1}(w_\ell)$. The neighbours of $v$ at level $s$ are its $b_{s-1}$ children, so
		\begin{equation}\label{eq:P-up}
			(P\Phi_s)(v)
			=\frac{b_{s-1}a_\ell}{\sqrt{m_{j,s}d_{s-1}d_s}}
			=\beta_{s-1}\Phi_{s-1}(v).
		\end{equation}
		
		\noindent(iii) Let $s=j+1$. The only possible nonzero entry at level $j$ is at $w$, and it equals
		\begin{equation}\label{eq:P-root}
			(P\Phi_{j+1})(w)
			=\frac{1}{\sqrt{d_jd_{j+1}}}\sum_{\ell=1}^{b_j}a_\ell=0.
		\end{equation}
		
		Since $(P\Phi_s)(v)=0$ for all other $v\in V$, \eqref{eq:P-down}--\eqref{eq:P-root} give, for $j+1<n$,
		\[
		\begin{aligned}
			P\Phi_{j+1}
			&=\beta_{j+1}\Phi_{j+2},\\
			P\Phi_s
			&=\beta_{s-1}\Phi_{s-1}+\beta_s\Phi_{s+1}\quad\text{for}\quad j+1<s<n,\\
			P\Phi_n
			&=\beta_{n-1}\Phi_{n-1}.
		\end{aligned}
		\]
		If $j+1=n$, then $P\Phi_n=0$ and $J_n=(0)$. Thus $\mathcal W(w,a)=\operatorname{span}\{\Phi_{j+1},\ldots,\Phi_n\}$ is invariant under $P$, and the restriction of $P$ to $\mathcal W(w,a)$ is represented by $J_{j+1}$ with respect to the ordered basis $(\Phi_{j+1},\ldots,\Phi_n)$.

		Since $\Psi_s$ is constant on $C_s$ and
	\begin{equation}\label{eq:phi-sum}
		\sum_{v\in C_s}\Phi_s(w,a)(v)
		=\sqrt{m_{j,s}}\sum_{\ell=1}^{b_j}a_\ell=0,
	\end{equation}
	we have $\langle\Phi_s(w,a),\Psi_s\rangle=0$. For $t\neq s$, $\Phi_s(w,a)(v)\,\Psi_t(v)=0$ for all $v\in V$. Hence  $\langle\Phi_s(w,a),\Psi_t\rangle=0$. Therefore $\mathcal W(w,a)\perp\mathcal R$. For the same $w$ and $p\neq q$, \eqref{eq:phi-inner} gives $\mathcal W(w,a^{(p)})\perp\mathcal W(w,a^{(q)})$.
	
	Now let $w\in C_j$ and $w'\neq w$. For $s\neq t$, $\langle\Phi_s(w,a),\Phi_t(w',a')\rangle=0$, so it suffices to take the same level $s$. If $D_s(w)\cap D_s(w')=\emptyset$, then $\Phi_s(w,a)(v)\,\Phi_s(w',a')(v)=0$ for all $v\in V$. Otherwise, say $D_s(w')\subseteq D_s(w_\ell)$. Then $\Phi_s(w,a)=a_\ell/\sqrt{m_{j,s}}$ on $D_s(w')$, and $\Phi_s(w',a')=0$ outside $D_s(w')$. Hence by \eqref{eq:phi-sum} for $(w',a')$,
	\[
	\langle\Phi_s(w,a),\Phi_s(w',a')\rangle
	=\frac{a_\ell}{\sqrt{m_{j,s}}}\sum_{v\in C_s}\Phi_s(w',a')(v)=0.
	\]
	This proves the mutual orthogonality, and the lemma follows.
	\end{proof}

		Two matrices $M$ and $N$ are orthogonally similar if $Q^{t}MQ=N$ for some orthogonal matrix $Q$. The notation $M\oplus N$ denotes the block diagonal matrix with diagonal blocks $M$ and $N$, and $M^{\oplus k}$ denotes the block diagonal matrix with $k$ copies of $M$. 
	\begin{proposition}\label{decomposition}
		Let $P$ be the discriminant matrix of the generalized Bethe tree $B(b_0,\ldots,b_{n-1})$. Then $P$ is orthogonally similar to
		\begin{equation}
			J_0\oplus
			\bigoplus_{\substack{0\leq j<n\\b_j\geq2}}
			J_{j+1}^{\,\oplus c_j(b_j-1)}.
		\end{equation}
	\end{proposition} 
	\begin{proof}
		For $0\leq s\leq n$, let $\mathcal B_s:=\{\Psi_s\}\cup\{\Phi_s(w,a^{(q)}) : 0\leq j<s,\ b_j\geq2,\ w\in C_j,\ 1\leq q\leq b_j-1\}$. 	Since $|C_j|=c_j$, each $j<s$ with $b_j\geq2$ contributes $c_j(b_j-1)$ vectors. Terms with $b_j=1$ are zero and using $c_{j+1}=c_jb_j$,
		\[
		|\mathcal B_s|=1+\sum_{j=0}^{s-1}c_j(b_j-1)
		=1+\sum_{j=0}^{s-1}(c_{j+1}-c_j)=c_s.
		\]
		Let $\mathcal B=\bigcup_{s=0}^n\mathcal B_s$. Then $|\mathcal B|=|V|$. The vectors $\Psi_0,\ldots,\Psi_n$ are orthonormal, so by Lemma~\ref{invariant-orthogonal}, $\mathcal B$ is an orthonormal set. Hence $\mathcal B$ is a real orthonormal basis of $\Cl^{|V|}$.
		
		Let $Q$ be the matrix with the vectors of $\mathcal B$ as columns, grouped by $\mathcal R$ and the subspaces $\mathcal W(w,a^{(q)})$. Then $Q$ is orthogonal. By Lemma~\ref{invariant-orthogonal}, each of these subspaces is $P$-invariant, and $P$ acts on $\mathcal R$ as $J_0$ and on $\mathcal W(w,a^{(q)})$ as $J_{j+1}$ for $w\in C_j$. For each $j$ with $b_j\geq2$, there are $c_j(b_j-1)$ such subspaces. So $Q^{t}PQ$ has the stated block diagonal form.
	\end{proof}
	\begin{lemma}\label{lem:block-support}
		Let $u\in C_s$, and let $J$ be $J_0$ or $J_{j+1}$ with $j<s$ and $b_j\geq2$. If $J$ has a $\mu$-eigenvector with nonzero coordinate at level $s$, then $\mu\in\Theta_P(u)$.
	\end{lemma}
	\begin{proof}
		An eigenvector $(y_0,\ldots,y_n)^{t}$ of $J_0$
		extends to the eigenvector $\sum_t y_t\Psi_t$ of $P$.
		Its entry at $u$ is $y_s/\sqrt{c_s}$.
		
		For $J_{j+1}$, let $w\in C_j$ be the ancestor of $u$, and
		let $w_\ell$ be its child whose subtree contains $u$.
		Choose one of the vectors $a^{(q)}$ from
		Proposition~\ref{decomposition} with $a^{(q)}_\ell\neq0$.
		Extending an eigenvector of $J_{j+1}$ to
		$\mathcal W(w,a^{(q)})$ gives entry
		$y_s a^{(q)}_\ell/\sqrt{m_{j,s}}$ at $u$.
		This is nonzero when $y_s\neq0$.
		For a symmetric matrix, an eigenvector with a nonzero entry
		at $u$ has eigenvalue in $\Theta_P(u)$.
	\end{proof}

	\subsection{Jacobi matrices}
	
	A \emph{Jacobi matrix} is a real symmetric tridiagonal matrix with positive entries immediately above and below the diagonal. By Proposition~\ref{decomposition}, the normalized adjacency matrix $P$ of a generalized Bethe tree is orthogonally similar to a direct sum of Jacobi matrices. We study some properties of Jacobi matrices which will be used subsequent sections.

	Let $J$ be a Jacobi matrix indexed by $0,\ldots,h$, and write $a_t=J_{t,t+1}$ for $0\leq t<h$. For $0\leq r\leq t\leq h$, let $J[r,t]$ denote the principal submatrix of $J$ indexed by $\{r,\ldots,t\}$. For completeness, we give a short proof of the following result (see also \cite[p.~268]{jacobi}).
	\begin{lemma}\label{jacobi-basic}
		Each eigenvalue of a Jacobi matrix $J$ is simple, and every eigenvector has nonzero first and last entries. If $v=(v_0,\ldots,v_h)^{t}$ is a $\lambda$-eigenvector of $J$ with $v_s=0$ for some $1\leq s\leq h-1$, then $(v_0,\ldots,v_{s-1})^{t}$ and $(v_{s+1},\ldots,v_h)^{t}$ are nonzero $\lambda$-eigenvectors of $J[0,s-1]$ and $J[s+1,h]$, respectively.
	\end{lemma}
	\begin{proof}
		Let $v=(v_0,\ldots,v_h)^{t}$ be a $\lambda$-eigenvector of $J$. Since the off-diagonal entries of $J$ are positive, the equation $Jv=\lambda v$ gives $v_j\in\operatorname{Span}\{v_0\}$ and $v_j\in\operatorname{Span}\{v_h\}$ for $0\leq j\leq h$. Thus $v_0\neq0$, $v_h\neq0$, and each eigenspace is one-dimensional. Hence each eigenvalue is simple.

		If $v_s=0$, restricting $Jv=\lambda v$ to the indices $0,\ldots,s-1$ and $s+1,\ldots,h$ gives the eigenvector equations for $J[0,s-1]$ and $J[s+1,h]$, respectively.
	\end{proof}
	For the rest of this section, we assume that $J$ has zero diagonal.  For $0\leq r\leq t\leq h$, we define $q_{r,t}(x):=\det(xI-2J[r,t])$. We also set $q_{r,r-1}(x)=1$ for $0\leq r\leq h+1$, using the convention that the determinant of the empty matrix is $1$. Then expansion along the last row gives
	\begin{equation}\label{eq:continuant}
		q_{r,t}(x)=xq_{r,t-1}(x)-4a_{t-1}^2q_{r,t-2}(x)\quad\text{for}\quad 0\leq r<t\leq h.
	\end{equation}
	Thus if squared off-diagonal entries of $J$ are rational, then $q_{r,t}\in\Ql[x]$. This holds for the blocks of $P$, since
	$\beta_i^2=b_i/(d_id_{i+1})$.

	A \emph{dyadic rational} is a number of the form $m/2^k$, where $m\in\Zl$ and $k$ is a nonnegative integer. We denote the ring of all such numbers by $\halfZ$. A complex number is \emph{integral over} a ring $R$ if it is a root of a monic polynomial in $R[x]$. A complex number integral over $\Zl$ is called an \emph{algebraic integer}. A subring $R$ of $\mathbb Q$ is \emph{integrally closed in $\mathbb Q$} if every rational number satisfying a monic polynomial with coefficients in $R$ belongs to $R$. We use the following standard coefficient criterion for $\Zl$ and $\halfZ$, both of which are integrally closed in $\Ql$.
	
	\begin{lemma}\label{dyadic-factor}
		Let $R=\Zl$ or $R=\halfZ$. If every root of a monic polynomial $f\in\Ql[x]$ is integral over $R$, then $f\in R[x]$. In particular, every monic factor in $\Ql[x]$ of a monic polynomial in $R[x]$ belongs to $R[x]$.
	\end{lemma}
	\begin{proof}
		Write $f(x)=(x-\alpha_1)\cdots(x-\alpha_m)=x^m+r_{m-1}x^{m-1}+\cdots+r_1x+r_0$. Each $r_i$ is a sum of products of the roots $\alpha_j$, up to sign, and is therefore integral over $R$. Since $a_i\in\Ql$ and $R$ is integrally closed in $\Ql$, we have $r_i\in R$, and hence $f\in R[x]$.
		
		Now let $g\in\Ql[x]$ be a monic factor of a monic polynomial $f\in R[x]$. Every root of $g$ is a root of $f$ and is therefore integral over $R$. Thus $g\in R[x]$ by the first assertion.
	\end{proof}
   
	\begin{lemma}\label{unsupported}
    Let $a_t^2\in\Ql$ for $0\leq t<h$, and let $0\leq s\leq h$. If $q_{s+1,h}\in\halfZ[x]$, and each eigenvalue of $2J$ that has an eigenvector with nonzero $s$-th coordinate is integral over $\halfZ$, then $q_{0,h}\in\halfZ[x]$.
	\end{lemma}
	\begin{proof}
		The characteristic polynomial $q_{0,h}$ of $2J$ is monic and has rational coefficients by \eqref{eq:continuant}. Hence, by Lemma~\ref{dyadic-factor}, it suffices to show that every eigenvalue of $2J$ is integral over $\halfZ$. This holds by assumption for every eigenvalue supported at $s$. If an eigenvalue $\lambda$ of $2J$ is unsupported at $s$, then Lemma~\ref{jacobi-basic}, applied to $2J$, gives $0<s<h$ and $q_{s+1,h}(\lambda)=0$.  Since $q_{s+1,h}$ is
		monic and belongs to $\halfZ[x]$, $\lambda$ is integral over
		$\halfZ$. This proves the result.
	\end{proof}
	Consider $J$ as the adjacency matrix of the path $0,\ldots,h$, with weight $a_t$ on the edge $\{t,t+1\}$. The entry $(J^k)_{ii}$ is the sum of the weights of all closed walks of length $k$ from $i$. Note that
	\begin{equation}\label{walk_2}
		(J^2)_{tt}=
		\begin{cases}
			a_0^2 & \text{if }t=0,\\
			a_{t-1}^2+a_t^2 & \text{if } 1\leq t\leq h-1,\\
			a_{h-1}^2 & \text{if } t=h.
		\end{cases}
	\end{equation}
	Two indices $i$ and $j$ of $J$ are \emph{cospectral} if $(J^k)_{ii}=(J^k)_{jj}$ for each integer $k\geq0$.
	\begin{lemma}\label{jacobi-symmetry}
		Let $i$ and $j$ be cospectral indices of $J$ with $i<j$. If $a_i=\cdots=a_{j-1}$, then $i+j=h$ and $a_t=a_{h-1-t}$ for $0\leq t<h$.  Moreover, $J$ is unchanged when its coordinates are taken in the revers order $h,\ldots,0$.
	\end{lemma}
	\begin{proof}
		The paths from $i$ to $0$ and from $j$ to $h$ have lengths $i$ and $h-j$, respectively. We show that these lengths are equal.
		
		Since $(J^2)_{ii}=(J^2)_{jj}$, $a_i=a_{j-1}$, and each edge weights are positive, \eqref{walk_2} gives $i=0$ if and only if $j=h$. If $i=0$, then $j=h$ and each edge weights are equal by hypothesis, and so the result follows. Otherwise, $0<i<j<h$, and \eqref{walk_2} gives $a_{i-1}^2+a_i^2=a_{j-1}^2+a_j^2$. Hence $a_{i-1}=a_j$.
		
		Let $r$ be the largest integer with $1\leq r\leq\min\{i,h-j\}$ such that $a_{i-\ell}=a_{j+\ell-1}$ for every $1\leq\ell\leq r$. Put $I_r=\{i-r,\ldots,j+r\}$. By the choice of $r$ and the hypothesis $a_i=\cdots=a_{j-1}$, we have $a_t=a_{i+j-1-t}$ for $i-r\leq t<j+r$. For each closed walk in $I_r$ starting at $i$, replace every index $k$ by $i+j-k$. This gives a closed walk starting at $j$ of the same length and weight. Applying the same rule again gives back the original walk.
		
		Consider closed walks of length $2r+2$. From $i$, the only possible walk that leaves $I_r$ takes $r+1$ steps to the left and return. From $j$, it takes $r+1$ steps to the right and returns. The walks that stay in $I_r$ contribute equally to $(J^{2r+2})_{ii}$ and $(J^{2r+2})_{jj}$. Hence, if exactly one of the two walks leaving $I_r$ exists, its positive weight contradicts cospectrality. If both exist, then their weights must agree, and hence
		\[
		\prod_{\ell=1}^{r+1}a_{i-\ell}^{\,2}
		=\prod_{\ell=1}^{r+1}a_{j+\ell-1}^{\,2}.
		\]
		Cancelling the first $r$ positive factors gives $a_{i-r-1}=a_{j+r}$, contradicting the maximality of $r$.
		
		Therefore neither walk exists, and hence $i-r=0$ and $j+r=h$. Thus $i+j=h$ and $I_r=\{0,\ldots,h\}$. The weight identity on $I_r$ now gives $a_t=a_{h-1-t}$ for $0\leq t<h$. Also, reversing the order of its rows and columns gives $J$ again.
	\end{proof}

\section{Restrictions on the transferring vertices}\label{sec:transfer-vertices}
 
The \emph{distance} $\dist(a,b)$ between two vertices $a$ and $b$ of a connected graph $G$ is the number of edges in a shortest path from $a$ to $b$.
\begin{lemma}\label{lem:automorphism-distances}
If a graph $G$ exhibits perfect state transfer from $u$ to $v$, then for each automorphism $\pi$ of $G$, $\dist(u,\pi(u))=\dist(v,\pi(v))$.
\end{lemma}
\begin{proof}
	Suppose $G$ exhibits perfect state transfer from $u$ to $v$ at time $\tau$. By Lemma~\ref{pre_lemma2}, $T_\tau(P)\eu=\ev$ and $T_\tau(P)\ev=\eu$. Since the permutation matrix of $\pi$ commutes with $P$, and so with $T_\tau(P)$. Therefore $T_\tau(P)\ee_{\pi(u)}=\ee_{\pi(v)}$. For each integer $\ell\geq0$, we have
	\[
	(P^\ell)_{\pi(u),u}=\langle\ee_{\pi(u)},P^\ell T_\tau(P)\ev\rangle=\langle T_\tau(P)\ee_{\pi(u)},P^\ell\ev\rangle=(P^\ell)_{\pi(v),v}.
	\]
	The entry $(P^\ell)_{a,b}$ is positive if and only if there is a walk of length $\ell$ from $a$ to $b$. This implies that  $\dist(a,b)=\min\{\ell\geq0:(P^\ell)_{a,b}>0\}$. Since $(P^\ell)_{\pi(u),u}=(P^\ell)_{\pi(v),v}$ for each $\ell\geq0$, we have the result.
\end{proof}
Let $w$ be a vertex of a rooted tree $G$. A \emph{child subtree} of $w$ is a component of $G-w$ containing a child of $w$, that is, a child of $w$ together with all its descendants. In a generalized Bethe tree, any two child subtrees of the same vertex are isomorphic as rooted trees.
\begin{lemma}\label{prop:different-levels}
	If a generalized Bethe tree exhibits perfect state transfer between
	vertices at different levels, then its underlying graph is a path.
\end{lemma}

\begin{proof}
	Let $G=B(b_0,\ldots,b_{n-1})$, and suppose $G$ exhibits perfect state transfer from $u\in C_i$ to $v\in C_j$, where $i<j$. We first show that $b_0=\cdots=b_{j-1}=1$.	Suppose, there exists $r$ is the first level with $r<j$ and $b_r\geq2$. Let $w$  is the unique vertex at level $r$.

	Let $W_1$ be the child subtree of $w$ containing $v$. If $u$ lies in another child subtree of $w$, choose it as $W_2$. Otherwise, choose any child subtree $W_2\neq W_1$. Then $W_1$ and $W_2$ are isomorphic. Choose a rooted isomorphism $\varphi:W_1\to W_2$. Let $\pi$ act as $\varphi$ on $W_1$, as $\varphi^{-1}$ on $W_2$, and fix all other vertices. Then $\pi$ is an automorphism of $G$. This implies that
	\[
	\dist(v,\pi(v))=2(j-r)\quad\text{and}\quad
	\dist(u,\pi(u))=
	\begin{cases}
		0      & i\leq r,\\
		2(i-r) & i>r,
	\end{cases}
	\]
	which contradicts Lemma~\ref{lem:automorphism-distances}, and hence $b_0=\cdots=b_{j-1}=1$.
	
	If $j=n$, then the graph $G$ is a path. Suppose that $j<n$. Then $\deg(v)\geq2$. By Lemma~\ref{pre_lemma2}, $\deg(u)=\deg(v)$, and hence $i>0$. Now $b_i=1$, so $d_i=d_j=2$ and $b_j=1$. Hence
	$\beta_i=\cdots=\beta_{j-1}=1/2$. Since $C_i=\{u\}$ and $C_j=\{v\}$, $\Psi_i=\eu$ and $\Psi_j=\ev$. Since $J_0$ represents
	$P|_{\mathcal R}$ and $u,v$ are strongly cospectral for $P$, we have $(J_0^\ell)_{i,i}=(P^\ell)_{u,u}
	=(P^\ell)_{v,v}=(J_0^\ell)_{j,j}$ for $\ell\geq0$. Therefore $i$ and $j$ are cospectral for $J_0$. Since the weights between $i$ and $j$ are equal, Lemma~\ref{jacobi-symmetry} shows that $J_0$ is invariant under reversal.
	
	Finally, let $q$ be the vector with $q_x=\sqrt{\deg(x)}$ for each
	$x\in V(G)$. Then $Pq=q$, since
	$(Pq)_x=\sum_{y\sim x}1/\sqrt{\deg(x)}=\sqrt{\deg(x)}$.
	With respect to the ordered basis $(\Psi_0,\ldots,\Psi_n)$ of
	$\mathcal R$, the $t$-th coordinate of $q$ is
	$z_t=\langle\Psi_t,q\rangle=\sqrt{c_t d_t}$. Hence $J_0z=z$. Reversing $z$ gives another positive $1$-eigenvector
	of the same norm. By Lemma~\ref{jacobi-basic}, the $1$-eigenspace
	is one-dimensional, so the two vectors are equal. In particular,
	$z_0=z_n$. Since $c_0=d_0=d_n=1$, we obtain
	\[
	\prod_{t=0}^{n-1}b_t=c_n=z_n^2=z_0^2=1.
	\]
	Every $b_t$ is a positive integer. Hence every $b_t=1$, and $G$
	is a path.
\end{proof}
 
We call a vertex $w$ an \emph{ancestor} of $v$ if $w$ lies on the
unique path from the root to $v$. In particular, each vertex is an
ancestor of itself. In our notation, for $w\in C_r$ and $v\in C_s$
with $r\leq s$, this is equivalent to $v\in D_s(w)$.
The \emph{lowest common ancestor} of $u$ and $v$ is the last vertex
shared by the paths from the root to $u$ and from the root to $v$.
 
\begin{lemma}\label{binary-split}
	Let $u,v$ be distinct vertices of $B$ in $C_s$ with lowest common ancestor $w\in C_r$. If $B$ exhibits perfect state transfer from $u$ to $v$ at time $\tau$, then $b_r=2$, $b_i=1$ for $r<i<s$, and $T_\tau(J_{r+1})\ee_s=-\ee_s$.
\end{lemma}
\begin{proof}
    Suppose $b_r\geq 3$. Let $W_1$ be the child subtree of $w$ containing $v$, and $W_2$ a child subtree of $w$ containing neither $u$ nor $v$. Since $W_1\cong W_2$, fix an isomorphism $\varphi:W_1\to W_2$. Define $\pi$ to be $\varphi$ on $W_1$, $\varphi^{-1}$ on $W_2$, and the identity elsewhere. Then $\pi$ is an automorphism of $B$ that fixes $u$ but not $v$, contradicting Lemma~\ref{pre_lemma2}(iv). Hence $b_r=2$.

    If $b_j\geq 2$ for some $r<j<s$, the same argument applied to the ancestor of $v$ at level $j$ gives an automorphism that fixes $u$ but not $v$, contradicting Lemma~\ref{pre_lemma2}(iv). Hence $b_i=1$ for $r<i<s$.
	
   Take $a=(1,-1)/\sqrt2$ at $w$. Since there is no branching between levels $r$ and $s$, \[\Phi_s(w,a)=\pm\frac{1}{\sqrt{2}}(\eu-\ev).\]
   By Lemma~\ref{pre_lemma2}(ii), $T_\tau(P)\eu=\ev$ and $T_\tau(P)\ev=\eu$, so $T_\tau(P)\Phi_s(w,a)=-\Phi_s(w,a)$. By Lemma~\ref{invariant-orthogonal}, $P$ acts on $\mathcal W(w,a)$ as $J_{r+1}$. Thus $T_\tau(J_{r+1})\ee_s=-\ee_s$.
\end{proof}

\section{Restrictions on the branching numbers}\label{sec:branching}

If perfect state transfer occurs between vertices at different levels, then the graph is a path, and so we restrict to transferring vertices at the same level. Throughout this section, we assume perfect state transfer occurs between $u$ and $v$ in $C_s$ with lowest common ancestor $w\in C_r$, where $r<s$. 

We first determine the branching below $w$, and then the path from the root to $w$. Both steps need the eigenvalues of $J_0$ and $J_{r+1}$ in polynomial form, which we set up now. Define polynomials $f_n,f_{n-1},\ldots,f_0$ by $f_n(x)=1$, $f_{n-1}(x)=x$, and
\begin{equation}\label{level-recurrence}
	f_{i-1}(x)=(b_i+1)x f_i(x)-b_i f_{i+1}(x) \quad\text{for}\quad1\leq i<n.
\end{equation}
Kubota et al.~\cite[Section~2.3]{bethetree} index these polynomials starting from the leaves. In our notation, their polynomials satisfy $g_k=f_{n-k}$ for $0\leq k\leq n$ and $g_{n+1}(x)=b_0\bigl(x f_0(x)-f_1(x)\bigr)$. Since $J_a$ is the submatrix of $J_0$ on levels $a,\ldots,n$, \cite[Lemma~2.5]{bethetree} identifies its eigenvalues with the zeros of $g_{n-a+1}$ for $0\leq a\leq n$. We state this result in our notation. For completeness, we include a short alternative proof.
\begin{lemma}\label{lem:level-recurrence}
	For $1\leq a\leq n$, $\lambda$ is an eigenvalue of $J_a$ if and only if $f_{a-1}(\lambda)=0$. Moreover, $\lambda$ is an eigenvalue of $J_0$ if and only if $\lambda f_0(\lambda)-f_1(\lambda)=0$.
\end{lemma}
\begin{proof}
	Let $y$ be a $\lambda$-eigenvector of $J_a$, where $0\leq a\leq n$. Write $y_i=\sqrt{c_i d_i}\,x_i$ for $a\leq i\leq n$. By Lemma~\ref{jacobi-basic}, $y_n\neq0$, so we may take $x_n=1$. For $a\geq1$, set $x_{a-1}=0$, since $J_a$ has no coordinate at level $a-1$.
	
	Using $c_{i+1}=c_i b_i$, we obtain $x_{n-1}=\lambda x_n$ and  $(b_i+1)\lambda x_i=x_{i-1}+b_i x_{i+1}$ for $\max\{1,a\}\leq i<n$.
	Starting from $x_n=1$ and working backwards, we obtain
	$x_i=f_i(\lambda)$ by \eqref{level-recurrence}.
	Thus, for $a\geq1$, the condition $x_{a-1}=0$ gives
	$f_{a-1}(\lambda)=0$. For $a=0$, the root equation
	$x_1=\lambda x_0$ gives $f_1(\lambda)=\lambda f_0(\lambda)$.
	
	Conversely, if the corresponding condition holds, then the vector with
	entries $y_i=\sqrt{c_i d_i}\,f_i(\lambda)$ for $a\leq i\leq n$
	satisfies all the eigenvalue equations. Since $f_n(\lambda)=1$,
	this vector is nonzero.
\end{proof} 
 
 \begin{lemma}\label{lower-shapes}
 	Exactly one of the following holds.
 	\begin{enumerate}[label=\textup{(\roman*)}]
 		\item $b_j=1$ for $r<j<n$.
 		\item For some positive integer $k$, $n=r+2k$ and $s\leq r+k$. The only branching level strictly between $r$ and $n$ is $r+k$, where $b_{r+k}=3$. Moreover, $\det(xI-2J_{r+1})=2T_{2k}(x/2)-1$.
 	\end{enumerate}
 \end{lemma}
 
 \begin{proof}
 	Let $\chi(x)=\det(xI-2J_{r+1})$. We first show that $\chi\in\Zl[x]$. Let $\mu$ be an eigenvalue of $2J_{r+1}$ with eigenvector $y:=(y_{r+1},\ldots,y_n)^{t}$, where the coordinates are indexed by levels.
 	
 	Suppose first that $y_s\neq0$. Then, by Lemma~\ref{binary-split}, $T_\tau(\mu/2)\,y_s=\ee_s^{t}\,T_\tau(J_{r+1})\,y=\bigl(T_\tau(J_{r+1})\ee_s\bigr)^{t}y=-y_s$. Thus $T_\tau(\mu/2)=-1$, so $\mu$ is a root of the monic polynomial $2T_\tau(x/2)+2\in\Zl[x]$. Hence $\mu$ is an algebraic integer.
 	
 	Now suppose that $y_s=0$. By Lemma~\ref{jacobi-basic} applied to $2J_{r+1}$, we have $s>r+1$, and $\mu$ is an eigenvalue of the leading principal submatrix of $2J_{r+1}$ corresponding to levels $r+1,\ldots,s-1$. By Lemma~\ref{binary-split}, $b_j=1$ for $r<j<s$, so $2\beta_j=1$ for $r+1\leq j\leq s-2$. Hence this submatrix has integer entries, and therefore $\mu$ is an algebraic integer.
 	
 	Thus every root of $\chi$ is an algebraic integer. By \eqref{eq:continuant}, $\chi$ is monic and $\chi\in\Ql[x]$. Hence $\chi\in\Zl[x]$ by Lemma~\ref{dyadic-factor}. Now \eqref{det_J_i} gives
 	\begin{equation}\label{branching-product}
 		\chi(2)=\frac{2^{n-r}}{\prod_{j=r+1}^{n-1}(b_j+1)}\in\Zl_{>0}.
 	\end{equation}
 	Consequently, $b_\ell+1$ is a power of $2$ for every $r+1\leq\ell\leq n-1$. There are $n-r-1$ such factors, each at least $2$, and their product divides $2^{n-r}$. Hence all these factors are $2$, except possibly one equal to $4$. Thus either \textup{(i)} holds or there is exactly one branching level $t$ with $r<t<n$, and then $b_t=3$.
 	
 	Suppose $b_t=3$. Since there is no branching level strictly between $r$ and $s$, we have $s\leq t$. Put $p=t-r$ and $q=n-t$. We now show that $p=q$.
 	
 	Since $b_i=1$ for $t<i<n$, \eqref{level-recurrence} gives $f_{i-1}(x)=2x f_i(x)-f_{i+1}(x)$ for $t<i<n$. As $f_n(x)=T_0(x)$ and $f_{n-1}(x)=T_1(x)$, we get $f_{n-\ell}(x)=T_\ell(x)$ for $0\leq\ell\leq q$. In particular, $f_t(x)=T_q(x)$ and $f_{t+1}(x)=T_{q-1}(x)$. Since $b_t=3$,
 	\[
 		f_{t-1}(x)=4xT_q(x)-3T_{q-1}(x)=2\bigl(2xT_q(x)-T_{q-1}(x)\bigr)-T_{q-1}(x)=2T_{q+1}(x)-T_{q-1}(x).
 	\]
 	With the convention $T_{-j}(x)=T_j(x)$ for $j\geq1$, we have
 	\[
 	f_{t-\ell}(x)=2T_{q+\ell}(x)-T_{q-\ell}(x)\quad\text{for}\quad0\leq\ell\leq p.
 	\]
 	Indeed, the two sides agree at $\ell=0,1$. Since $b_{t-\ell}=1$ for $1\leq\ell<p$, both sides satisfy the same Chebyshev recurrence in $\ell$, and the identity follows by induction. For $\ell=p$, this gives $f_r(x)=2T_{p+q}(x)-T_{|p-q|}(x)$.

 	Note that $2T_m(x/2)$ is monic of degree $m$ in $\Zl[x]$ for $m\geq1$. Since $|p-q|<p+q$, the polynomial $f_r(x/2)=2T_{p+q}(x/2)-T_{|p-q|}(x/2)$ is monic of degree $p+q$. By Lemma~\ref{lem:level-recurrence}, $f_r(x/2)=0$ if and only if $x/2$ is an eigenvalue of $J_{r+1}$, that is, $x$ is an eigenvalue of $2J_{r+1}$. Since $J_{r+1}$ has order $n-r=p+q$ and its eigenvalues are distinct by Lemma~\ref{jacobi-basic}, $f_r(x/2)$ and $\chi(x)$ are monic with the same $p+q$ distinct roots. Hence
 	\begin{equation}\label{eq:unequal-tail-polynomial}
 		\chi(x)=2T_{p+q}(x/2)-T_{|p-q|}(x/2).
 	\end{equation}
 	Since $\chi(x)$ and $2T_{p+q}(x/2)$ belong to $\Zl[x]$, \eqref{eq:unequal-tail-polynomial} gives $T_{|p-q|}(x/2)\in\Zl[x]$. For $m\geq1$, $T_m(x/2)$ has leading coefficient $1/2$, so this implies $p=q=k$ for some positive integer $k$. Hence $t=r+k$, $n=r+2k$, and $s\leq r+k$. Finally, since $T_0=1$, \eqref{eq:unequal-tail-polynomial} gives $\chi(x)=2T_{2k}(x/2)-1$, which completes the proof.
 \end{proof}

\begin{lemma}\label{no-earlier-branch}
	There is no branching before level $r$, that is, $b_0=\cdots=b_{r-1}=1$.
\end{lemma}

\begin{proof}
	For $s+1\leq i<n$, Lemma~\ref{lower-shapes} gives $b_i\in\{1,3\}$, so $d_i\in\{2,4\}$. Also $d_n=1$. Hence $d_i\in\{1,2,4\}$ for $s+1\leq i\leq n$, and so $\beta_i^2=b_i/(d_id_{i+1})\in\halfZ$ for $s+1\leq i<n$. By \eqref{eq:continuant}, it follows that $\det(xI-2J_{s+1})\in\halfZ[x]$.	If $s=n$, the empty matrix has characteristic polynomial $1$.
	
	Suppose $b_j\geq2$ for some $j<r$.  Let $\mu$ be an eigenvalue of $J_{j+1}$ with an eigenvector whose coordinate at level $s$ is nonzero. By Lemma~\ref{lem:block-support}, $\mu\in\Theta_P(u)$, and hence $2\mu$ is an algebraic integer by Lemma~\ref{pre_lemma2}(v). Now Lemma~\ref{unsupported}, applied to $J_{j+1}$, gives $\det(xI-2J_{j+1})\in\halfZ[x]$. Hence $\det(2I-2J_{j+1})\in\halfZ$.
	
	On the other hand, $J_{j+1}$ has order $n-j$, so \eqref{det_J_i} gives
	\[
	\det(2I-2J_{j+1})=2^{n-j}\det(I-J_{j+1})=\frac{2^{n-j}}{\kappa},\quad\text{where } \kappa:=\prod_{\ell=j+1}^{n-1}(b_\ell+1).
	\]
	Since $j<r<n$ and $b_r=2$ by Lemma~\ref{binary-split}, $\kappa$ contains the factor $b_r+1=3$. If $2^{n-j}/\kappa=m/2^k$ with $m\in\Zl$ and $k\geq0$, then $2^{n-j+k}=m\kappa$, so $3$ divides a power of $2$, a contradiction. Hence $b_j=1$ for each $j<r$.
\end{proof}

Let $F_m(x)$ be defined by $F_0(x)=0$, $F_1(x)=1$, and $F_{m+1}(x)=xF_m(x)-F_{m-1}(x)$ for $m\geq1$. Set $F_{-m}(x)=-F_m(x)$ for $m\geq1$. For $m\geq1$, the polynomial $F_m$ is monic of degree $m-1$ with integer coefficients. The recurrence also gives
\begin{equation}\label{eq:sine-polynomials}
	F_m(2\cos\theta)=\frac{\sin(m\theta)}{\sin\theta}
	\quad\text{for }\sin\theta\neq0.
\end{equation}
For every positive integer $k$ and every integer $m$,
\begin{equation}\label{crb1844}
	F_{k+m}(x)-F_{k-m}(x)=2T_k(x/2)F_m(x).
\end{equation}
 
By Lemma~\ref{no-earlier-branch}, there is no branching before level $r$. By Lemma~\ref{binary-split}, we have $b_r=2$. Lemma~\ref{lower-shapes} gives two possibilities for the branching numbers at levels $r+1,\ldots,n-1$.  If $r=0$ in case~\textup{(i)} of Lemma~\ref{lower-shapes}, then $b_0=2$ and $b_j=1$ for $1\leq j<n$, so the generalized Bethe tree is a path. The other three cases are those of the next lemma, and for them we compute the characteristic polynomial of $2J_0$.

The formulas in the following lemma are special cases of the characteristic polynomial calculations of Kubota et al.~\cite[Section~3]{bethetree}, they used the substitution $x=(z+z^{-1})/2$. We derive these formulas directly from the level recurrence \eqref{level-recurrence} and express $\det(xI-2J_0)$ in terms of $F_m$.

\begin{lemma}\label{lem:radial-polynomials}
	Let $b_0=\cdots=b_{r-1}=1$ and $b_r=2$. Then following are hold.
	\begin{enumerate}[label=\textup{(\roman*)}]
	\item If $r\geq1$, $n=r+k$, and $b_j=1$ for $r<j<n$, then
	\begin{equation}\label{eq:radial-two-arms}
		\det(xI-2J_0)
		=(x^2-4)\left(F_{r+k}(x)+\frac13 F_{k-r}(x)\right).
	\end{equation}
	\item If $r\geq1$, $n=r+2k$, $b_{r+k}=3$, and $b_j=1$ for
	$r<j<n$ with $j\neq r+k$, then
	\begin{equation}\label{eq:radial-branched-arms}
		\det(xI-2J_0)
		=(x^2-4)\left(
		F_{r+2k}(x)+\frac23 T_k(x/2)F_{k-r}(x)
		\right).
	\end{equation}
	\item If $r=0$, $n=2k$, $b_k=3$, and $b_j=1$ for
	$0<j<n$ with $j\neq k$, then
	\begin{equation}\label{eq:radial-root-split}
		\det(xI-2J_0)=(x^2-4)F_{2k}(x).
	\end{equation}
\end{enumerate}
\end{lemma}
\begin{proof}
	Let $D(x)=\det(xI-2J_0)$, $\lambda=x/2$, and $g(\lambda)=\lambda f_0(\lambda)-f_1(\lambda)$, where the polynomials $f_i$ are given by \eqref{level-recurrence}. Note that $\deg f_i=n-i$ and the leading coefficient of $f_{i-1}$ is $b_i+1$ times that of $f_i$ for $1\leq i<n$. Thus $g$ has degree $n+1$ and leading coefficient $\prod_{i=1}^{n-1}(b_i+1)$. By Lemmas~\ref{lem:level-recurrence} and~\ref{jacobi-basic},  roots of $g$ are the $n+1$ distinct eigenvalues of $J_0$. Hence $\det(\lambda I-J_0)=g(\lambda)/\prod_{i=1}^{n-1}(b_i+1)$. Therefore
	\begin{equation}\label{crb1843}
		D(x)=2^{n+1}\det(\lambda I-J_0)
		=4g(\lambda)\prod_{i=1}^{n-1}\frac{2}{b_i+1}.
	\end{equation}
	
	We also use the identity
	\begin{equation}\label{crb1837}
		\lambda T_j(\lambda)-T_{j-1}(\lambda)
		=(\lambda^2-1)F_j(x)
		=\frac{x^2-4}{4}F_j(x)
		\quad\text{for } j\geq1.
	\end{equation}
	Indeed, for $\lambda=\cos\theta$ with $0<\theta<\pi$, all three expressions in \eqref{crb1837} equal $-\sin\theta\sin(j\theta)$.

	Suppose that $r\geq1$. Since $b_i=1$ for $1\leq i<r$, we have
	\[
	\begin{aligned}
		T_i(\lambda)f_{i-1}(\lambda)-T_{i-1}(\lambda)f_i(\lambda)
		&=\bigl(2\lambda T_i(\lambda)-T_{i-1}(\lambda)\bigr)f_i(\lambda)-T_i(\lambda)f_{i+1}(\lambda)\\
		&=T_{i+1}(\lambda)f_i(\lambda)-T_i(\lambda)f_{i+1}(\lambda).
	\end{aligned}
	\]
	Thus $T_r(\lambda)f_{r-1}(\lambda)-T_{r-1}(\lambda)f_r(\lambda)=T_1(\lambda)f_0(\lambda)-T_0(\lambda)f_1(\lambda)=g(\lambda)$.
	Since $b_r=2$, we have $f_{r-1}=3\lambda f_r-2f_{r+1}$, and hence
	\begin{equation}\label{crb1845}
		\begin{aligned}
			g(\lambda)
			&=T_r(\lambda)\bigl(3\lambda f_r(\lambda)-2f_{r+1}(\lambda)\bigr)-T_{r-1}(\lambda)f_r(\lambda)\\
			&=\bigl(\lambda T_r(\lambda)-T_{r-1}(\lambda)\bigr)f_r(\lambda)
			+2T_r(\lambda)\bigl(\lambda f_r(\lambda)-f_{r+1}(\lambda)\bigr).
		\end{aligned}
	\end{equation}
	
	\begin{enumerate}[label=\textup{(\roman*)}]
		\item Since $b_i=1$ for $r<i<n$ and $n=r+k$, the level recurrence, with $f_n=T_0(\lambda)$ and $f_{n-1}=T_1(\lambda)$, gives $f_r=T_k(\lambda)$ and $f_{r+1}=T_{k-1}(\lambda)$. Here \eqref{crb1843} gives $D(x)=\frac{8}{3}g(\lambda)$. Substituting these values in \eqref{crb1845} and using \eqref{crb1837}, we obtain
		\[
		\begin{aligned}
			D(x)
			&=\frac{2(x^2-4)}{3}
			\bigl(T_k(\lambda)F_r(x)+2T_r(\lambda)F_k(x)\bigr)\\
			&=(x^2-4)\left(F_{r+k}(x)+\frac{F_{k-r}(x)}{3}\right).
		\end{aligned}
		\]
		The last equality follows from \eqref{crb1844}, since $2T_k(\lambda)F_r(x)=F_{r+k}(x)-F_{k-r}(x)$ and $2T_r(\lambda)F_k(x)=F_{r+k}(x)+F_{k-r}(x)$. This proves \eqref{eq:radial-two-arms}.
		
		\item By the calculation in the proof of Lemma~\ref{lower-shapes}, $f_r=2T_{2k}(\lambda)-1$ and $f_{r+1}=2T_{2k-1}(\lambda)-T_1(\lambda)$.	Since $T_1(\lambda)=\lambda$, we have $\lambda f_r-f_{r+1}=2\bigl(\lambda T_{2k}(\lambda)-T_{2k-1}(\lambda)\bigr)$. Since $r\geq1$, the factors at levels $r$ and $r+k$ in \eqref{crb1843} are $2/3$ and $1/2$, respectively. Thus $D(x)=\frac{4}{3}g(\lambda)$. Using \eqref{crb1845} and \eqref{crb1837}, we obtain
		\[
		D(x)=\frac{x^2-4}{3}
		\Bigl(\bigl(2T_{2k}(\lambda)-1\bigr)F_r(x)+4T_r(\lambda)F_{2k}(x)\Bigr).
		\]
		By \eqref{crb1844}, $2T_{2k}(\lambda)F_r(x)=F_{2k+r}(x)-F_{2k-r}(x)$ and $4T_r(\lambda)F_{2k}(x)=2\bigl(F_{r+2k}(x)+F_{2k-r}(x)\bigr)$. Hence
		\[
		D(x)=(x^2-4)\left(F_{r+2k}(x)+\frac{F_{2k-r}(x)-F_r(x)}{3}\right).
		\]
		Taking $m=k-r$ in \eqref{crb1844} gives $F_{2k-r}(x)-F_r(x)=2T_k(x/2)F_{k-r}(x)$. This proves \eqref{eq:radial-branched-arms}.
		
		\item Here $r=0$. The formulas for $f_r$ and $f_{r+1}$ in \textup{(ii)} also hold for $r=0$, and give $f_0=2T_{2k}(\lambda)-1$ and $f_1=2T_{2k-1}(\lambda)-T_1(\lambda)$. Hence, by \eqref{crb1837},
		\[
		g(\lambda)=\lambda f_0(\lambda)-f_1(\lambda)
		=2\bigl(\lambda T_{2k}(\lambda)-T_{2k-1}(\lambda)\bigr)
		=\frac{x^2-4}{2}F_{2k}(x).
		\]
		Since $b_k=3$ and $b_i=1$ for $1\leq i<n$ with $i\neq k$, \eqref{crb1843} gives $D(x)=2g(\lambda)$. Thus $D(x)=(x^2-4)F_{2k}(x)$, proving \eqref{eq:radial-root-split}.\qedhere
	\end{enumerate}
\end{proof}

\begin{proposition}\label{necessary-families}
	If a generalized Bethe tree exhibits perfect state transfer, then its underlying graph is isomorphic to a path, $S_k(B(2,3))$, or $S_k(B(1,2,3))$ for some positive integer $k$.
\end{proposition}
\begin{proof}
	Suppose that $B(b_0,\ldots,b_{n-1})$ exhibits perfect state transfer from $u$ to $v$. If $u$ and $v$ lie at different levels, then by Lemma~\ref{prop:different-levels} the graph is a path. Thus we may assume that $u,v\in C_s$, and let $w\in C_r$ be their lowest common ancestor. By Lemmas~\ref{binary-split} and~\ref{no-earlier-branch}, we have $b_r=2$ and there is no branching before level $r$. If $r=0$, the two cases of Lemma~\ref{lower-shapes} give a path and $S_k(B(2,3))$, respectively. Hence we may assume that $r\geq1$.
	
	Let $\mu$ be an eigenvalue of $J_0$ with an eigenvector whose coordinate at level $s$ is nonzero. By Lemma~\ref{lem:block-support}, $\mu\in\Theta_P(u)$, so $2\mu$ is an algebraic integer by Lemma~\ref{pre_lemma2}(v). Moreover, as shown in the proof of Lemma~\ref{no-earlier-branch}, $\det(xI-2J_{s+1})\in\halfZ[x]$. Thus Lemma~\ref{unsupported}, applied to $J_0$, gives $D(x):=\det(xI-2J_0)\in\halfZ[x]$.
	
	By \eqref{eq:radial-two-arms} and \eqref{eq:radial-branched-arms}, we have
	\[
	D(x)-(x^2-4)F_n(x)=
	\begin{cases}
		\frac13 (x^2-4)F_{k-r}(x) & \text{in case \textup{(i)} of Lemma~\ref{lower-shapes}},\\
		\frac23 (x^2-4)F_{k-r}(x)T_k(x/2) & \text{in case \textup{(ii)} of Lemma~\ref{lower-shapes}}.
	\end{cases}
	\]
	Since $F_n\in\Zl[x]$ and $D(x)\in\halfZ[x]$, $D(x)-(x^2-4)F_n(x)\in\halfZ[x]$. If $r\neq k$, then $F_{k-r}$ has leading coefficient $\pm1$. Since $2T_k(x/2)$ is monic, the leading coefficient of $D(x)-(x^2-4)F_n(x)$ is  $\pm1/3$ in either case, a contradiction. Hence $r=k$.

	In case~\textup{(i)} of Lemma~\ref{lower-shapes}, the tree consists of three paths of length $k$ meeting at $w$. For $1\leq i\leq k$, let $x_i$ and $y_i$ be the vertices at distance $i$ from $w$ on the path containing $v$ and on the path from $w$ to the root, respectively. Consider the permutation $\pi=(x_1\,y_1)\cdots(x_k\,y_k)$.  Then $\pi$ is an automorphism of the underlying graph that fixes $u$ and sends $v=x_{s-r}$ to $y_{s-r}\neq v$, contrary to Lemma~\ref{pre_lemma2}(iv).

	In case~\textup{(ii)} of Lemma~\ref{lower-shapes}, we have $n=3k$, and the only branching levels are $k$ and $2k$, with branching numbers $2$ and $3$, respectively. Thus the underlying graph is $S_k(B(1,2,3))$.
\end{proof}

\section{Perfect state transfer on generalized Bethe trees}\label{sec:transfer-times}

We now prove that each graph in Proposition~\ref{necessary-families} exhibits perfect state transfer, and determine all transferring pairs and their minimum transfer times. Label the  vertices of $P_{n+1}$ as $v_0,\ldots,v_n$ in path order.
\begin{proposition}\label{prop:path-transfers}
   The path $P_{n+1}$ exhibits perfect state transfer occurs between $v_i$ and $v_j$ if and only if $i+j=n$. The minimum transfer time for each such pair is $n$.
\end{proposition}

\begin{proof}
	Let $\delta_j=\deg(v_j)$. For $0\leq\ell\leq n$, put $\theta_\ell=\ell\pi/n$ and $\mu_\ell=\cos\theta_\ell$. Define $y^{(\ell)}$ by $y^{(\ell)}_j=\sqrt{\delta_j}\cos(j\theta_\ell)$ for $0\leq j\leq n$. Then one can show that
	\[
	(Py^{(\ell)})_j=\frac{1}{\sqrt{\delta_j}}\sum_{v_i\sim v_j}\cos(i\theta_\ell)=\mu_\ell y^{(\ell)}_j.
	\]
	Thus $Py^{(\ell)}=\mu_\ell y^{(\ell)}$. Since the $\mu_\ell$ are distinct, the vectors $y^{(0)},\ldots,y^{(n)}$  form an eigenbasis of $P$.
	
	Let $Q$ be the permutation matrix defined by $Q\ee_{v_j}=\ee_{v_{n-j}}$. Since $\delta_{n-j}=\delta_j$ and $\cos((n-j)\theta_\ell)=(-1)^\ell\cos(j\theta_\ell)$, we have $Qy^{(\ell)}=(-1)^\ell y^{(\ell)}$. Also, $T_n(\mu_\ell)=\cos(\ell\pi)=(-1)^\ell$. Thus $T_n(P)=Q$. By Lemma~\ref{pre_lemma1}, perfect state transfer occurs  between distinct vertices $v_j$ and $v_{n-j}$ at time $n$.
	
	Suppose that perfect state transfer occurs between $v_i$ and $v_j$ at time $\tau$, where $i\neq j$. By Lemma~\ref{pre_lemma2}, $v_i$ and $v_j$ are strongly cospectral and $\deg v_i=\deg v_j$. Since $\mu_1$ is simple, $y^{(1)}_i=\pm y^{(1)}_j$, and hence $\cos(i\pi/n)=\pm\cos(j\pi/n)$.  Since cosine is strictly decreasing on $[0,\pi]$ and $i\neq j$, we must have $\cos(i\pi/n)=-\cos(j\pi/n)$. Hence $i+j=n$. In particular, $y^{(1)}_i=-y^{(1)}_j\neq0$. Since $T_\tau(P)$ is symmetric, taking the inner product of both sides of $T_\tau(P)\ee_{v_i}=\ee_{v_j}$ with $y^{(1)}$ gives $T_\tau(\mu_1)y_i^{(1)}=y_j^{(1)}$. Thus $\cos(\tau\pi/n)=-1$, so $\tau$ is an odd multiple of $n$. Hence the minimum transfer time is $n$.
\end{proof}

We now consider $S_k(B(2,3))$ and $S_k(B(1,2,3))$, where $k\geq1$. In each tree, let $w\in C_r$ be the unique vertex with two children. Then $r=0$ in $S_k(B(2,3))$ and $r=k$ in $S_k(B(1,2,3))$. In both trees, $n=r+2k$, $b_{r+k}=3$, and $b_i=1$ for $r<i<n$ with $i\neq r+k$. Recall that the rows and columns of $J_{r+1}$ are indexed by the levels $r+1,\ldots,n$.

\begin{lemma}\label{lem:contrast-time}
	For $S_k(B(2,3))$ or $S_k(B(1,2,3))$, $T_{6k}(J_{r+1})=-I$. Moreover, if $T_\tau(J_{r+1})\ee_s=-\ee_s$ for some $r+1\leq s\leq n$ and some positive integer $\tau$, then $\tau$ is an odd multiple of $6k$.
\end{lemma}

\begin{proof}
	For either tree, the calculation leading to \eqref{eq:unequal-tail-polynomial} gives
	\begin{equation}\label{eq:contrast-polynomial}
		\det(xI-2J_{r+1})=2T_{2k}(x/2)-1.
	\end{equation}
	Since $J_{r+1}$ is a block of $P$, its eigenvalues lie in $[-1,1]$. For $0\leq\theta\leq\pi$, it follows from \eqref{eq:contrast-polynomial} that $\mu=\cos\theta$ is an eigenvalue of $J_{r+1}$ if and only if $\cos(2k\theta)=1/2$. Hence, for each eigenvalue $\mu$ of $J_{r+1}$, we have $T_{6k}(\mu)=T_3(T_{2k}(\mu))=T_3(1/2)=-1$. Since $J_{r+1}$ is symmetric, $T_{6k}(J_{r+1})=-I$.
	
	The smallest such $\theta$ is $\pi/(6k)$, so the largest eigenvalue of $J_{r+1}$ is $\rho=\cos(\pi/(6k))$. Since $J_{r+1}$ is nonnegative and irreducible, the Perron--Frobenius theorem gives a positive $\rho$-eigenvector $y$. Taking the inner product of both sides of $T_\tau(J_{r+1})\ee_s=-\ee_s$ with $y$ gives $T_\tau(\rho)y_s=-y_s$. Since $y_s>0$, we have $\cos(\tau\pi/(6k))=-1$. Hence $\tau$ is an odd multiple of $6k$.
\end{proof}

For both $S_k(B(2,3))$ and $S_k(B(1,2,3))$, Lemmas~\ref{prop:different-levels} and~\ref{binary-split} show that any two transferring vertices lie at the same level $s$ and have $w\in C_r$ as their lowest common ancestor. Since $b_i=1$ for $r<i<s$ and $b_{r+k}=3$, we have $r<s\leq r+k$.

For such $s$, take $a=(1,-1)/\sqrt{2}$. We may label the two vertices in $C_s$ as $u$ and $v$ so that
\[
\Psi_s=\frac{\eu+\ev}{\sqrt{2}}\quad\text{and}\quad\Phi_s(w,a)=\frac{\eu-\ev}{\sqrt{2}}.
\]
By Lemmas~\ref{pre_lemma1} and~\ref{pre_lemma2}(i), perfect state transfer occurs in these trees between $u$ and $v$ at time $\tau$ if and only if $T_\tau(P)\Psi_s=\Psi_s$ and $T_\tau(P)\Phi_s(w,a)=-\Phi_s(w,a)$. By \eqref{eq:radial-matrix} and Lemma~\ref{invariant-orthogonal}, perfect state transfer occurs between $u$ and $v$ at time $\tau$ if and only if
\begin{equation}\label{crb0032}
	T_\tau(J_0)\ee_s=\ee_s
	\quad\text{and}\quad
	T_\tau(J_{r+1})\ee_s=-\ee_s.
\end{equation}

\begin{proposition}\label{prop:b23-transfers}
	In $S_k(B(2,3))$, perfect state transfer occurs precisely between the two vertices of degree $4$ at the minimum time $6k$.
\end{proposition}
\begin{proof}
	In $S_k(B(2,3))$, $r=0$ and $n=2k$. Hence by Lemma~\ref{lower-shapes}, $1\leq s\leq k$. 	By \eqref{eq:radial-root-split}, we have $\det(xI-2J_0)=(x^2-4)F_{2k}(x)$. Hence, by \eqref{eq:sine-polynomials}, the eigenvalues of $J_0$ are $\mu_\ell=\cos(\ell\pi/2k)$ for $0\leq\ell\leq2k$. By \eqref{crb0032} and Lemma~\ref{lem:contrast-time}, any transfer time has the form $\tau=(2m+1)6k$ for some $m\geq0$. Therefore,
	\begin{equation}\label{eq:b23-radial-sign}
		T_\tau(\mu_\ell)=\cos\bigl(3(2m+1)\ell\pi\bigr)=(-1)^\ell.
	\end{equation}
	
	Let $y$ be a $\mu_\ell$-eigenvector of $J_0$, and write $y_i=\sqrt{c_i d_i}\,z_i$. Since $y_0\neq0$ by Lemma~\ref{jacobi-basic}, we may take $z_0=1$. As in the proof of Lemma~\ref{lem:level-recurrence}, the root equation gives $z_1=\mu_\ell$. Since $b_i=1$ for $1\leq i<k$, we have $z_{i+1}=2\mu_\ell z_i-z_{i-1}$. Thus $z_s=T_s(\mu_\ell)=\cos(s\ell\pi/2k)$ for $0\leq s\leq k$. 
	
	By \eqref{eq:b23-radial-sign} and the eigenvector formula,	$T_\tau(J_0)\ee_s=\ee_s$ if and only if $\cos(s\ell\pi/(2k))=0$ for every odd $\ell$. Taking $\ell=1$ forces $s=k$, since $1\leq s\leq k$. Hence, by \eqref{crb0032}, perfect state transfer can occur only at level $k$.
	
    At $s=k$, we have $z_k=\cos(\ell\pi/2)=0$ for every odd $\ell$. Thus $\ee_k$ is orthogonal to every $\mu_\ell$-eigenvector of $J_0$ with odd $\ell$. Since $T_{6k}(\mu_\ell)=1$ for every even $\ell$ by \eqref{eq:b23-radial-sign}, it follows that $T_{6k}(J_0)\ee_k=\ee_k$. Also, Lemma~\ref{lem:contrast-time} gives $T_{6k}(J_1)\ee_k=-\ee_k$. Hence by \eqref{crb0032}, $S_k(B(2,3))$ exhibits perfect state transfer between the two vertices in $C_k$  at time $6k$. Lemma~\ref{lem:contrast-time} gives the minimum time $6k$.
\end{proof}
 
In $S_k(B(1,2,3))$, let $c$ be the unique vertex of degree $3$ and let $x,y$ be the two vertices of degree $4$. For $0\leq i\leq k$, let $x_i$ and $y_i$ denote the vertices at distance $i$ from $c$ on the paths from $c$ to $x$ and $y$, respectively.

\begin{proposition}\label{prop:b123-transfers}
	In $S_k(B(1,2,3))$, perfect state transfer occurs precisely between $x_i$ and $y_i$ for $1\leq i\leq k$ at the minimum time $6k$.
\end{proposition}

\begin{proof}
	Here $w=c$, $r=k$, and $n=3k$. For $1\leq j\leq k$, the two vertices in $C_{k+j}$ are $x_j$ and $y_j$. Thus these are the only possible pairs. Since $F_0=0$, putting $r=k$ in \eqref{eq:radial-branched-arms} gives $\det(xI-2J_0)=(x^2-4)F_{3k}(x)$. By \eqref{eq:sine-polynomials}, the eigenvalues of $J_0$ are $\cos(\ell\pi/(3k))$ for $0\leq\ell\leq3k$. Since $T_{6k}(\cos(\ell\pi/(3k)))=\cos(2\ell\pi)=1$, we have $T_{6k}(J_0)=I$. Also, $T_{6k}(J_{k+1})=-I$ by Lemma~\ref{lem:contrast-time}. Taking $s=k+j$ in \eqref{crb0032}, we obtain perfect state transfer between $x_j$ and $y_j$ at time $6k$ for every $1\leq j\leq k$. By Lemma~\ref{lem:contrast-time}, $6k$ is the minimum transfer time.
\end{proof}

\begin{proof}[\textbf{Proof of Theorem~\ref{pst_gbt}}]
	By Proposition~\ref{necessary-families}, any generalized Bethe tree that exhibits perfect state transfer belongs to one of the three families listed in the theorem. Conversely, Propositions~\ref{prop:path-transfers}, \ref{prop:b23-transfers}, and~\ref{prop:b123-transfers} show that every graph in these families exhibits perfect state transfer.
\end{proof}

\begin{proof}[\textbf{Proof of Proposition~\ref{prop:pst-pairs}}]
	Propositions~\ref{prop:path-transfers}, \ref{prop:b23-transfers}, and~\ref{prop:b123-transfers} give all transferring pairs and their minimum transfer times for $P_{n+1}$, $S_k(B(2,3))$, and $S_k(B(1,2,3))$, respectively, as listed in Table~\ref{tab:pst-pairs}.
\end{proof}

\end{document}